\documentclass[MSNbibl,number,citesort,seceqn,dvips]{arxbj}
\usepackage{dcolumn}
\usepackage{graphicx}

\aid{0}
\volume{19}
\issue{5A}
\pubyear{2013}
\firstpage{2120}
\lastpage{2151}
\doi{10.3150/12-BEJ446} 

\makeatletter
\newcommand{\rrvert}{\vert}
\newcommand{\llvert}{\vert}
\newproclaim{assumption}{Assumption}[section]
\newtheorem{theorem}{Theorem}[section]
\newremark{remark}{Remark}[section]
\newtheorem{lemma}{Lemma}[section]

\newcommand{\vech}{\operatorname{vech}}
\newcommand{\var}{\operatorname{var}}
\newcommand{\diag}{\operatorname{diag}}

\newcommand{\inP}{\stackrel{{P}}{\longrightarrow}}
\newcommand{\inD}{\stackrel{\mathcal{D}}{\longrightarrow}}
\newcommand{\as}{\stackrel{\mathrm{a.s.}}{\longrightarrow}}
\newcommand{\Int}{\int_0^1}

\newcommand{\eqref}[1]{(\ref{#1})}
\newcolumntype{d}[1]{D{.}{.}{#1}}
\makeatother

\begin{document}
\begin{frontmatter}

\title{A test of significance in functional quadratic regression}
\runtitle{Functional quadratic regression}

\begin{aug}
\author{\fnms{Lajos} \snm{Horv\'ath}} \and
\author{\fnms{Ron} \snm{Reeder}\corref{}\ead[label=e2]{reeder@math.utah.edu}}
\runauthor{L. Horv\'ath and R. Reeder} 
\address{Department of Mathematics, University of Utah, Salt Lake City, UT 84112, USA.\\\printead{e2}}
\end{aug}

\received{\smonth{4} \syear{2011}}
\revised{\smonth{4} \syear{2012}}

%
\begin{abstract}
We consider a quadratic functional regression model in which a scalar
response depends on a functional predictor; the common functional
linear model is a special case. We wish to test the significance of the
nonlinear term in the model. We develop a testing method which is based
on projecting the observations onto a suitably chosen finite
dimensional space using functional principal component analysis. The
asymptotic behavior of our testing procedure is established. A
simulation study shows that the testing procedure has good size and
power with finite sample sizes. We then apply our test to a data set
provided by Tecator, which consists of near-infrared absorbance spectra
and fat content of meat.
\end{abstract}

%

\end{frontmatter}

\section{Introduction and results}\label{sintro5}

In a predictive model, it may be more natural and appropriate for
certain quantities to be represented as trajectories rather than a
single number (Kirkpatrick and Heckman~\cite{kirkpatrickheckman1989}). For example, a young
animal's size may be considered as a function of time, giving a growth
trajectory. A~model to predict a certain response from growth
trajectories is useful to animal breeders because they may be able to
produce more valuable animals by changing their growth patterns (Fitzhugh \cite{fitzhugh1976}).
M\"uller and Zhang \cite{muellerzhang2005} used egg-laying
trajectories from Mediterranean fruit flies to predict a female fly's
remaining lifetime. Frank and Friedman \cite{frankfriedman1993} and
Wold \cite{wold1993} provide an early discussion on the applications of principal
components to analyze curves in chemistry. Further examples for
analysis of data when the observations are curves can be found in
Fan and Lin \cite{fanlin1998}, Laukaitis and Ra\v{c}kauskas \cite{lauka2005},
Cardot \textit{et al.} \cite{cardotprchalsarda2007} and Zhang and Chen \cite{zhangchen2007}. For surveys on
functional data analysis, we refer to the books of Ferraty and Vieu \cite{ferratyvieu2006}, Ferraty and Romain \cite{ferratyromain2011} and
Horv\'ath and Kokoszka \cite{horvathkokoszka2012}.

\begin{figure}

\includegraphics{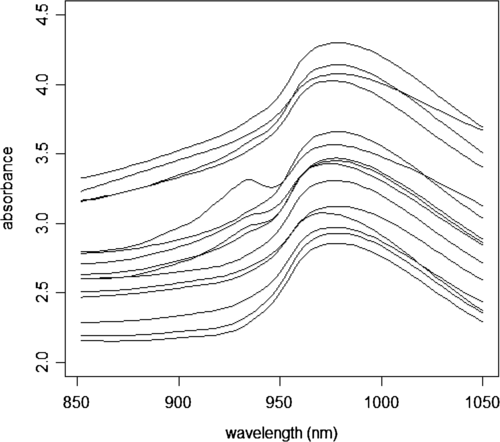}

\caption{Absorbance trajectories from 15
samples of finely chopped pure meat.}\label{fspectra}
\end{figure}

Yao and M\"uller \cite{muelleryao2010} and Borggaard and Thodberg \cite{borggaardthodberg1992} used
absorbance trajectories to predict the fat content of meat samples. The
absorbance at any particular wavelength is a measurement related to the
proportion of light that passes through a meat sample. A~representative
sample of 15 of the 240 absorbance trajectories are pictured in
Figure~\ref{fspectra}.

In functional regression, special attention has been given to
functional linear models (Cardot \textit{et al.} \cite{cardotfms2003}, Shen and Faraway \cite{shenfaraway2004}, Cai and Hall \cite{caihall2006},
Hall and Horowitz \cite{hallhorowitz2007}).
However, it is pointed out in Yao and M\"uller \cite{muelleryao2010} that this
model imposes a constraint on the regression relationship that may not
be appropriate in some scenarios. Yao and M\"uller \cite{muelleryao2010}
generalized this to a functional polynomial model, which has greater
flexibility. %
In functional polynomial regression, as in standard polynomial
regression, one must balance the costs and benefits of using more
parameters in the model. In this paper, we will develop a test to
determine if a quadratic term is justified in the model or if a
functional linear model adequately describes the regression relationship.

We assume that the predictor, $X_n(t)$, is defined on a finite interval
which, without loss of generality, will be $[0,1]$. %
The functional quadratic model in which a scalar response, $Y_n$, is
paired with a functional predictor, $X_n(t)$, is defined as
%
\begin{equation}\label{eqmodel}
Y_n = \mu+ \Int k(t) X_n^c(t)\,\mathrm{d}t+
\Int\Int h(s,t) X_n^c(s) X_n^c(t)
\,\mathrm{d}t\,\mathrm{d}s+\varepsilon_n,
\end{equation}
where $X_n^c(t)= X_n(t)- E (X_n(t) )$ is the centered
predictor process and $\varepsilon_n$ is a random error. The functions
$k(t)$ and $h(s,t)$ are regression parameter functions in the model
\eqref{eqmodel}. If $h(s,t)=0$, then $\mu= E(Y_n)$ and \eqref{eqmodel} reduces to the functional linear model
%
\begin{equation}
\label{eqlin} Y_n = \mu+ \Int k(t) X_n^c(t)\,\mathrm{d}t+
\varepsilon_n.
\end{equation}

Yao and M\"uller \cite{muelleryao2010} developed estimators for the functions $k$
and $h$ and prediction theory for the model \eqref{eqmodel}. Cardot and Sarda \cite{cardotsarda2011} and Mas and Pumo \cite{maspumo2011}
point out in their
survey papers that since we can choose a function in \eqref{eqlin},
the functional linear model can be used in a large variety of
applications. The functional linear model provides a very simple
relation between $X_n(t)$ and $Y_n$, so it is important to check if the
more involved quadratic model \eqref{eqmodel} provides a real
improvement. In other words, one should test whether the quadratic term
is really needed. To test the significance of the quadratic term in
\eqref{eqmodel}, we test the null hypothesis,
%
\begin{equation}
\label{eqnull5} H_0\dvt  h(s,t)=0,
\end{equation}
against the alternative
\[
H_A\dvt h(s,t)\neq0.
\]
To reduce the dimensionality and avoid overfitting in our functional
regression model, we will project the predictor process onto a suitably
chosen finite dimensional space. The space is spanned by the
eigenfunctions of $C(t,s)=E(X_n(t)-\mu_X(t))(X_n(s)-\mu_X(s))$, the
covariance function of the predictor process, where $\mu_X(t)=EX_n(t)$.
We will denote the eigenfunctions and associated eigenvalues by $\{
(v_i(t), \lambda_i), 1\le i \le\infty\}$. We can and will assume that
$\lambda_i$ is the $i$th largest eigenvalue and that the
eigenfunctions are orthonormal.
It is clear that we can assume that $h$ is symmetric, and we also
impose the condition that the kernels are in $L^2$:
\begin{eqnarray}
\label{eqh} h(s,t)&=&h(t,s) \quad\mbox{and}\quad
 \Int\Int h^2(s,t)\,\mathrm{d}t\,\mathrm{d}s<\infty,\\
\label{eqk} \Int k^2(t)\,\mathrm{d}t&<&\infty.
\end{eqnarray}
Thus, we have the expansions
%
\begin{eqnarray}
\label{eqexph} %
h(s,t)&=&\sum_{i=1}^{\infty}
\sum_{j=1}^{\infty} a_{i,j}
v_j(s) v_i(t)
\nonumber
\\[-8pt]
\\[-8pt]
&=&\sum_{i=1}^{\infty} a_{i,i}
v_i(s) v_i(t) + \sum_{i=1}^{\infty}
\sum_{j=i+1}^{\infty} a_{i,j}
\bigl(v_j(s) v_i(t) + v_i(s)
v_j(t) \bigr)
\nonumber
\end{eqnarray}
and
%
\begin{equation}
\label{eqexpk} k(t)=\sum_{i=1}^{\infty}
b_{i} v_i(t).
\end{equation}
Let $\langle\cdot, \cdot\rangle$ denote the inner product in $L^2$.
By projecting onto the space spanned by $\{v_1, \dots, v_p\}$ and using
\eqref{eqexph} and \eqref{eqexpk}, we can write the model \eqref
{eqmodel} as
%
\begin{equation}
\label{eqmodel1} %
Y_n = \mu+\sum
_{i=1}^{p} b_{i} \bigl\langle
X_n^c, v_i \bigr\rangle+ \sum
_{i=1}^{p} \sum_{j=i}^{p}
\bigl(2-1\{i=j\}\bigr) a_{i,j} \bigl\langle X_n^c,
v_i \bigr\rangle \bigl\langle X_n^c,
v_j \bigr\rangle + \varepsilon_n^*, %
\end{equation}
where
\begin{eqnarray*}
\varepsilon_n^*&=& \varepsilon_n + \sum
_{i=p+1}^{\infty} b_{i} \bigl\langle
X_n^c, v_i \bigr\rangle+\sum
_{i=p+1}^{\infty} \sum_{j=i}^{\infty}
\bigl(2-1\{i=j\}\bigr) a_{i,j} \bigl\langle X_n^c,
v_i \bigr\rangle \bigl\langle X_n^c,
v_j \bigr\rangle\\
&&{}+\sum_{i=1}^{p}
\sum_{j=p+1}^{\infty} 2 a_{i,j} \bigl
\langle X_n^c, v_i \bigr\rangle \bigl\langle
X_n^c, v_j \bigr\rangle.
\end{eqnarray*}
We note that \eqref{eqmodel1} is written as a standard linear model,
but the error term, $\varepsilon^*_n$, and the design points, $\{
\langle
X_n^c, v_i\rangle, 1\leq i \leq p\}$, are dependent.

We observe $(Y_n, X_n(t), 0\leq t \leq1), 1\leq n \leq N.$ %
Unfortunately, we cannot use \eqref{eqmodel1} directly for
statistical inference since $v_i(t)$ and $\mu_X(t)$ are unknown. We
estimate $\mu_X(t)$ and $C(t,s)$ with the corresponding empiricals
\[
\bar{X}(t)=\frac{1}{N}\sum_{n=1}^N
X_n(t)
\]
and
\[
\hat{C}(t,s)=\frac{1}{N}\sum_{n=1}^N
\bigl(X_n(t)-\bar{X}(t) \bigr) \bigl(X_n(s)-\bar{X}(s)
\bigr).
\]
The eigenvalues and the corresponding eigenfunctions of $\hat{C}(t,s)$
are denoted by $\hat{\lambda}_1\geq\hat{\lambda}_2\geq\cdots$ and
$\hat{v}_1, \hat{v}_2,\dots$\,. Eigenfunctions corresponding to unique
eigenvalues are uniquely determined up to signs. For this reason, we
cannot expect more than to have $\hat{c}_{i} \hat{v}_{i}(t)$ be close
to $v_i(t)$, where the $\hat{c}_{i}$'s are random signs. %
We replace equation \eqref{eqmodel1} with
%
\begin{eqnarray}
\label{eqmodel2} Y_n &=& \mu+\sum_{i=1}^{p}
b_{i} \langle X_n - \bar{X}, \hat{c}_i
\hat{v}_i \rangle\nonumber\\[-8pt]\\[-8pt]
&&{}+ \sum_{i=1}^{p}
\sum_{j=i}^{p} \bigl(2-1\{i=j\}\bigr)
a_{i,j} \langle X_n - \bar{X}, \hat{c}_i
\hat{v}_i \rangle \langle X_n - \bar{X},
\hat{c}_j \hat {v}_j \rangle + \varepsilon_n^{**},\nonumber
\end{eqnarray}
where
\begin{eqnarray*}
\varepsilon_n^{**}&=&\varepsilon_n^* + \sum
_{i=1}^{p} b_{i} \bigl\langle
X_n^c, v_i-\hat{c}_i
\hat{v}_i \bigr\rangle+ \sum_{i=1}^p
b_i \langle\bar{X}-\mu_X,\hat{c}_i
\hat{v}_i \rangle
\\
&&{}- \sum_{i=1}^{p} \sum
_{j=i}^{p} \bigl(2-1\{i=j\}\bigr) a_{i,j}
\bigl( \langle X_n - \bar{X}, \hat{c}_i
\hat{v}_i \rangle \langle X_n - \bar{X}, \hat
{c}_j \hat {v}_j \rangle-\bigl\langle
X_n^c, v_i\bigr\rangle \bigl\langle
X_n^c, v_j\bigr\rangle \bigr). %
\end{eqnarray*}
We can write \eqref{eqmodel2} in the concise form\vspace*{-2pt}

\begin{equation}
\label{eqconcise} {\mathbf Y} = \hat{\mathbf Z} %
\pmatrix{ \tilde{
\mathbf A}
\cr
\tilde{\mathbf B}
\cr
\mu} %
+ {\boldsymbol
\varepsilon^{**}},
\end{equation}
where
\begin{eqnarray*}
{\mathbf Y}&=& %
(Y_1, Y_2, \dots,
Y_N) %
^T\in R^N,\\
 \tilde{\mathbf A}&=&
\vech \bigl( \bigl\{ \hat{c}_i \hat{c}_j
a_{i,j} \bigl(2-1\{ i=j\} \bigr), 1\le i \le j\le p \bigr
\}^T \bigr)\in R^{p(p+1)/2},
\\
\tilde{\mathbf B}&=& %
( \hat{c}_1 b_1,
\hat{c}_2 b_2, \dots, \hat{c}_p
b_p ) %
^T\in R^p,\qquad {\boldsymbol
\varepsilon}^{**}= %
\bigl( \varepsilon_1^{**},
\varepsilon_2^{**}, \dots, \varepsilon_N^{**}
\bigr) %
^T\in R^N,
\end{eqnarray*}
and $\hat{{\mathbf Z}}$ is a $N\times(p(p+1)/2 +p+1)$ matrix given by
\[
\hat{{\mathbf Z}}= %
\pmatrix{ \hat{\mathbf D}_1^T
& \hat{\mathbf F}_1^T &1
\cr
\hat{\mathbf
D}_2^T & \hat{\mathbf F}_2^T &1
\cr
\vdots& \vdots&\vdots
\vspace*{3pt}\cr
\hat{\mathbf D}_N^T & \hat{
\mathbf F}_N^T &1 } %
\]
with
\begin{eqnarray*}
\hat{\mathbf D}_n &=& \vech \bigl(\bigl\{ \langle
\hat{v}_i, X_n - \bar{X} \rangle \langle
\hat{v}_j , X_n - \bar{X} \rangle, 1 \le i \le j \le p
\bigr\}^T \bigr)\in R^{p(p+1)/2},
\\
\hat{\mathbf F}_n &=& %
\bigl( \langle X_n -
\bar{X}, \hat{v}_1 \rangle, \langle X_n - \bar{X},
\hat{v}_2 \rangle, \dots, \langle X_n - \bar{X},
\hat{v}_p \rangle \bigr) %
^T\in R^p.
\end{eqnarray*}
The half-vectorization, $\vech(\cdot)$, stacks the columns of the lower
triangular portion of the matrix under each other.
Although we write our model in the form of a general linear model, it
is important to note that it is not a classical linear model. First,
${\boldsymbol\varepsilon^{**}}$ is correlated with $\hat{\mathbf Z}$
because ${\boldsymbol\varepsilon^{**}}$ contains additional error
terms which come from projecting onto a p-dimensional space. Another
important difference between \eqref{eqconcise} and a classical linear
model is that the parameters to be estimated, $\tilde{\mathbf A}$ and
$\tilde{\mathbf B}$, are random; they depend on the random signs,
$\hat
{c}_{i}$. We estimate $\tilde{\mathbf A}$, $\tilde{\mathbf B}$, and
$\mu
$ using the least squares estimator:
%
\begin{equation}
\label{eqhats} %
\pmatrix{ \hat{\mathbf A}
\cr
\hat{\mathbf
B}
\cr
\hat{\mu} } %
= \bigl(\hat{\mathbf Z}^T\hat{\mathbf
Z} \bigr)^{-1}\hat{\mathbf Z}^T{\mathbf Y}. %
\end{equation}
To represent elements of $\hat{\mathbf A}$ and $\hat{\mathbf B}$, we
will use the notation that $\hat{\mathbf A} = \vech(\{\hat
{a}_{i,j}(2-1\{i=j\}), 1\le i \le j \le p\}^T)\in R^{p(p+1)/2}$ and
$\hat{\mathbf B} =
(\hat{b}_1,  \hat{b}_2,  \dots,  \hat{b}_{p})
^T\in R^p$.

We expect, under $H_0$, that $\hat{\mathbf A}$ will be close to zero
since $\tilde{\mathbf A}$ is zero. If $H_0$ is not correct, we expect
the magnitude of $\hat{\mathbf A}$ to be relatively large. This
suggests that a testing procedure could be based on $\hat{\mathbf A}$.
Due to the random signs coming from the estimation of the
eigenfunctions, $\hat{\mathbf A}$ will not be asymptotically normal.
However, if the random signs are ``taken out,'' asymptotic normality can
be established. Hence, our test statistic will be a quadratic form of
$\hat{\mathbf A}$ with some random weight matrices. Let%
\begin{eqnarray*}
\hat{\mathbf G} &=& \frac{1}{N} \sum_{n=1}^{N}
\hat{\mathbf D}_n \hat {\mathbf D}_n^T,
\\
\hat{\mathbf M} &=& \frac{1}{N} \sum_{n=1}^N
\hat{\mathbf D}_n,
\end{eqnarray*}
and
\[
\hat{\tau}^2 = \frac{1}{N}\sum_{n=1}^N
\hat{\varepsilon}_n^2,
\]
where
\[
\hat{\varepsilon}_n= Y_n -\hat{\mu} - \sum
_{i=1}^{p} \hat{b}_{i} \langle
X_n - \bar{X}, \hat{v}_i \rangle- \sum
_{i=1}^p \sum_{j=i}^p
\bigl(2-1\{ i=j\}\bigr) \hat{a}_{i,j} \langle X_n -
\bar{X}, \hat{v}_i \rangle \langle X_n - \bar{X}, \hat
{v}_j \rangle
\]
are the residuals under $H_0$. We reject the null hypothesis if
\[
U_N=\frac{N}{\hat{\tau}^2} \hat{\mathbf A}^T \bigl(\hat{
\mathbf G}- \hat {\mathbf M}\hat{\mathbf M}^T\bigr) \hat{\mathbf A}
\]
is large. The main result of this paper is the asymptotic distribution
of $U_N$ under the null hypothesis.
First, we discuss the assumptions needed to establish asymptotics for~$U_N$:
%
\begin{assumption}\label{agaussian}
$\{ X_n(t), n\ge1 \}$ is a sequence of independent, identically
distributed Gaussian processes.
\end{assumption}
%
\begin{assumption}\label{aXmoments}
\[
E \biggl( \Int X_n^2(t)\, \mathrm{d}t \biggr)^4 <
\infty.
\]
\end{assumption}
%
\begin{assumption}\label{avarepsilon}
$\{\varepsilon_n\}$ is a sequence of independent, identically
distributed random variables satisfying $E\varepsilon_n = 0$ and $E
\varepsilon_n^4 < \infty$,
\end{assumption}
\noindent and
%
\begin{assumption}\label{aindep} The sequences $\{\varepsilon_n\}$ and
$\{X_n(t)\}$ are independent.
\end{assumption}

The last condition is standard in functional data analysis. It implies
that the eigenfunctions $v_1, v_2, \ldots, v_p$ are unique up to a sign.

\begin{assumption}\label{aunique}
\[
\lambda_1 > \lambda_2 > \cdots> \lambda_{p+1}.
\]
\end{assumption}

\begin{theorem}\label{thmain}
If $H_0$, \eqref{eqk} and Assumptions \ref{agaussian}--\ref{aunique}
are satisfied, then
\[
U_N \inD\chi^2(r),
\]
where $r=p(p+1)/2$ is the dimension of the vector $\hat{\mathbf A}$.
\end{theorem}

The proof of Theorem \ref{thmain} is given in Section \ref{sproofs5}.
%
\begin{remark}\label{r1} %
By the Karhunen--Lo\`eve expansion, every centered, square integrable
process, $X_n^c(t)$, can be written as
\[
X^c_n(t)=\sum_{\ell=1}^\infty
\xi_{n,\ell} \varphi_\ell(t),
\]
where $\varphi_\ell$ are orthonormal functions. Assumption \ref
{agaussian} can be replaced with the requirement that $\xi_{n,1}$,
$\xi_{n,2}$, $\ldots$\,, $\xi_{n,p}$ are independent with $E\xi_{n,\ell}^3=0$
and $E\xi_{n,\ell}=0$ for all $1\leq\ell\leq p$.
\end{remark}

Our last result provides a simple condition for the consistency of the
test based on $U_N$. Let ${\mathbf A} = \vech(\{{a}_{i,j}(2-1\{i=j\}
), 1\le i \le j \le p\}^T)$, that is, the first $r=p(p+1)/2$
coefficients in the expansion of $h$ in \eqref{eqexph}.
%
\begin{theorem}\label{thconst}
If \eqref{eqh}, \eqref{eqk}, Assumptions \ref{agaussian}--\ref
{aunique} are satisfied and ${\mathbf A}\neq{\mathbf0},$ then we
have that
\[
U_N \inP\infty.
\]
\end{theorem}

The condition ${\mathbf A}\neq{\mathbf0}$ means that $h$ is not the 0
function in the space spanned by the functions $v_i(t)v_j(s), 1\leq i,
j\leq p$. Our alternative means that if the quadratic term is needed to
explain the $p$-dimensional
projections, then the test will see it.

\section{A simulation study}\label{sempirical5}

In this section, we investigate the empirical size and power of the
testing procedure for finite sample sizes. Seeking to obtain a test of
size $\alpha= 0.01$, $0.05$, or $0.10$, a rejection region was chosen
according to the limiting distribution of the test statistic. Since the
limiting distribution is $\chi^2(r)$, the rejection region is $(\Delta,
\infty)$, where $P(\chi^2(r) > \Delta) = \alpha$.
Simulated data was then used to compute the outcome of the test
statistic. Iterating this procedure 5000 times, we kept track of the
proportion of times that the outcome fell in the predetermined
rejection region. When simulations are done under $H_0$, this gives us
the empirical size of the test, which we expect to be close to the
nominal size, $\alpha$, for large sample sizes. When simulations are
done under the alternative, $H_A$, the proportion gives us the
empirical power of the test.

In our first simulation study, the $\varepsilon_n$'s were generated
according to the distribution of independent standard normals. We
generated the $X_n(t)$'s according to the distribution of independent
standard Brownian motions. Then, using $k(t) = 1$ and $h(s,t)=c$, we
obtained $Y_n$ according to \eqref{eqmodel}. Thus, the power of the
test is a function of the parameter $c$. In particular, when $c=0$, the
null hypothesis is true. The resulting empirical size and power are
given in Table \ref{tempirical}.

\begin{table}
\tablewidth=8cm
\tabcolsep=0pt
\caption[Power of test using iid observations]{Empirical power of test
(in \%) based on 5000 simulations using i.i.d. Brownian motions for
$X_n(t)$ and i.i.d. standard normals for $\varepsilon_n$ and $N=200$}\label{tempirical}
\begin{tabular*}{8cm}{@{\extracolsep{\fill}}llll@{}}
\hline
$c$& $p=1$ & $p=2$ & $p=3$\\ \hline 
& \multicolumn{3}{c}{$\alpha=0.01$} \\
$0.0$ &$ \hphantom{1}1.02 $&$ \hphantom{1}1.37 $&$ \hphantom{1}1.95$\\ 
$0.2$ &$ 10.81 $&$ \hphantom{1}6.87 $&$ \hphantom{1}6.52$\\ 
$0.4$ &$ 49.51 $&$ 37.24 $&$ 29.76$\\ 
$0.6$ &$ 86.68 $&$ 77.74 $&$ 70.19$\\ 
$0.8$ &$ 98.50 $&$ 96.05 $&$ 92.98$\\ 
$1.0$ &$ 99.94 $&$ 99.57 $&$ 99.05$\\ 
[5pt]
&  \multicolumn{3}{c}{$\alpha=0.05$} \\
$0.0$ &$ \hphantom{1}5.15 $&$ \hphantom{1}6.00 $&$ \hphantom{1}7.44 $\\ 
$0.2$ &$ 25.90 $&$ 19.17 $&$ 18.02 $\\ 
$0.4$ &$ 72.10 $&$ 60.31 $&$ 50.38 $\\ 
$0.6$ &$ 95.21 $&$ 90.43 $&$ 85.77 $\\ 
$0.8$ &$ 99.60 $&$ 98.90 $&$ 97.60 $\\ 
$1.0$ &$ 99.99 $&$ 99.87 $&$ 99.84 $\\ 
[5pt]
&  \multicolumn{3}{c}{$\alpha=0.10$} \\
$0.0$ &$ 10.27 $&$ 11.18 $&$ 13.35$\\ 
$0.2$ &$ 36.60 $&$ 29.50 $&$ 27.03 $\\ 
$0.4$ &$ 80.89 $&$ 71.08 $&$ 62.27 $\\ 
$0.6$ &$ 97.60 $&$ 94.77 $&$ 90.91 $\\ 
$0.8$ &$ 99.85 $&$ 99.47 $&$ 98.57 $\\ 
$1.0$ &$ 99.99 $&$ 99.95 $&$ 99.91 $\\ 
\hline
\end{tabular*}
\end{table}

The distribution of our test statistic has been shown to converge to a
$\chi^2(r)$. Thus, we expect the empirical and nominal size to be close
for samples of size $N=200$. 
Since our testing procedure depends on the choice of how many principal
components to keep, results are given in Table~\ref{tempirical} for
$p=1$, $2$, and $3$. One possible method of selecting $p$ is to follow
the advice of Ramsay and Silverman \cite{ramsaysilverman2005} and choose $p$ so that
approximately 85\% of the variance within a sample is described by the
first $p$ principal components.

\begin{table}
\tablewidth=8cm
\tabcolsep=0pt
\caption[Power of test using non-Gaussian observations]{Empirical power
of test (in \%) based on 5000 simulations using non-Gaussian $X_n(t)$
and non-normal $\varepsilon_n$ and $N=200$}\label{tviolations}
\begin{tabular*}{8cm}{@{\extracolsep{\fill}}llll@{}}
\hline
$c$& $p=1$ & $p=2$ & $p=3$\\ \hline
& \multicolumn{3}{c}{$\alpha=0.01$} \\
$0.0$ &$ \hphantom{11}2.40 $&$ \hphantom{11}1.20 $&$ \hphantom{11}1.85 $\\ 
$0.2$ &$ \hphantom{1}57.70 $&$ \hphantom{1}46.75 $&$ \hphantom{1}37.50 $\\ 
$0.4$ &$ \hphantom{1}96.90 $&$ \hphantom{1}95.55 $&$ \hphantom{1}91.20 $\\ 
$0.6$ &$ \hphantom{1}99.90 $&$ 100.00 $&$ \hphantom{1}99.70 $\\ 
$0.8$ &$ 100.00 $&$ 100.00 $&$ 100.00 $\\ 
$1.0$ &$ 100.00 $&$ 100.00 $&$ 100.00 $\\ 
[5pt]
& \multicolumn{3}{c}{$\alpha=0.05$} \\
$0.0$ &$ \hphantom{11}8.00 $&$ \hphantom{11}5.75 $&$ \hphantom{11}8.15 $\\
$0.2$ &$ \hphantom{1}74.50 $&$ \hphantom{1}64.55 $&$ \hphantom{1}56.45 $\\
$0.4$ &$ \hphantom{1}99.40 $&$ \hphantom{1}98.35 $&$ \hphantom{1}96.55 $\\
$0.6$ &$ \hphantom{1}99.95 $&$ 100.00 $&$ \hphantom{1}99.85 $\\
$0.8$ &$ 100.00 $&$ 100.00 $&$ 100.00 $\\
$1.0$ &$ 100.00 $&$ 100.00 $&$ 100.00 $\\
[5pt]
& \multicolumn{3}{c}{$\alpha=0.10$} \\
$0.0$ &$ \hphantom{1}13.60 $&$ \hphantom{1}12.15 $&$ \hphantom{1}14.60 $\\
$0.2$ &$ \hphantom{1}82.30 $&$ \hphantom{1}74.25 $&$ \hphantom{1}65.55 $\\
$0.4$ &$ \hphantom{1}99.65 $&$ \hphantom{1}99.10 $&$ \hphantom{1}97.95 $\\
$0.6$ &$ \hphantom{1}99.95 $&$ 100.00 $&$ \hphantom{1}99.90 $\\
$0.8$ &$ 100.00 $&$ 100.00 $&$ 100.00 $\\
$1.0$ &$ 100.00 $&$ 100.00 $&$ 100.00 $\\
\hline
\end{tabular*}
\end{table}

%
\begin{table}
\tablewidth=10cm
\tabcolsep=0pt
\caption{Empirical power of test (in \%) based on 5000 simulations
using i.i.d. Brownian motions for $X_n(t)$, i.i.d. standard normals for
$\varepsilon_n$, $h(s,t)=c$ and $N=500$}\label{t1}
\begin{tabular*}{10cm}{@{\extracolsep{\fill}}ld{3.2}d{3.2}d{3.2}d{3.2}d{3.2}@{}}
\hline
$c$&  \multicolumn{1}{l}{$p=1$} & \multicolumn{1}{l}{$p=2$} & \multicolumn{1}{l}{$p=3$} & \multicolumn{1}{l}{$p=4$} & \multicolumn{1}{l}{$p=5$}\\ \hline
& \multicolumn{5}{c}{$\alpha=0.01$} \\
$0.0$ &1.1 &  1.3  &  1.15  &  1.6  &  2  \\
$0.2$ &30.35 &  20.35  &  12.85  &  10.85  &  9.05  \\
$0.4$ &91.9 &  84.25  &  74.35  &  67.7  &  58.7  \\
$0.6$ &100&  99.7  &  98.75  &  98  &  96.75  \\
$0.8$ &100&  100  &  100  &  99.9  &  99.85  \\
$1.0$ &100&  100  &  100  &  100  &  100  \\[5pt]
& \multicolumn{5}{c}{$\alpha=0.05$} \\
$0.0$ & 5.6 &  5.75  &  6.05  &  6.4  &  8.05  \\
$0.2$ & 53.05 &  40  &  31.35  &  26.3  &  24.85  \\
$0.4$ & 97.9 &  93.7  &  88.55  &  85.85  &  78.2  \\
$0.6$ & 100 &  99.9  &  99.6  &  99.6  &  99.1  \\
$0.8$ & 100 &  100  &  100  &  100  &  99.95  \\
$1.0$ & 100 &  100  &  100  &  100  &  100  \\[5pt]
& \multicolumn{5}{c}{$\alpha=0.10$} \\
$0.0$ & 10.6 &  11.05  &  11.55  &  11.85  &  13.5  \\
$0.2$ & 65 &  52.45  &  43.75  &  37.65  &  35.85  \\
$0.4$ & 99.3 &  96.6  &  93.1  &  91.4  &  85.7  \\
$0.6$ & 100 &  99.95  &  99.75  &  99.8  &  99.75  \\
$0.8$ & 100 &  100  &  100  &  100  &  99.95  \\
$1.0$ & 100 &  100  &  100  &  100  &  100  \\ \hline
\end{tabular*}
\vspace*{6pt}
\end{table}


%
\begin{table}
\tablewidth=10cm
\tabcolsep=0pt
\caption{Empirical power of test (in \%) based on 5000 simulations
using non-Gaussian $X_n(t)$ and non-normal $\varepsilon_n$, $h(s,t)=c$
and $N=500$}\label{t2}
\begin{tabular*}{10cm}{@{\extracolsep{\fill}}ld{3.2}d{3.2}d{3.2}d{3.2}@{}}
\hline
$c$&  \multicolumn{1}{l}{$p=1$} & \multicolumn{1}{l}{$p=2$} & \multicolumn{1}{l}{$p=3$} & \multicolumn{1}{l}{$p=4$} \\ \hline
& \multicolumn{4}{c}{$\alpha=0.01$} \\
0.0  &  1.45  &  1.35  &  1.55  &  1.75  \\
0.2  &  90.3  &  82.55  &  14.1  &  77.65  \\
0.4  &  100  &  100  &  65.35  &  100  \\
0.6  &  100  &  100  &  96  &  100  \\
0.8  &  100  &  100  &  99.9  &  100  \\
1.0  &  100  &  100  &  100  &  100  \\[5pt]
& \multicolumn{4}{c}{$\alpha=0.05$} \\
0.0  &  5.45  &  6.1  &  7.05  &  6.6  \\
0.2  &  96.3  &  92  &  28.35  &  90.05  \\
0.4  &  100  &  100  &  81.55  &  100  \\
0.6  &  100  &  100  &  98.75  &  100  \\
0.8  &  100  &  100  &  100  &  100  \\
1.0  &  100  &  100  &  100  &  100  \\[5pt]
& \multicolumn{4}{c}{$\alpha=0.10$} \\
0.0  &  10.35  &  11.5  &  13.1  &  12.7  \\
0.2  &  97.7  &  95.25  &  39.2  &  93.8  \\
0.4  &  100  &  100  &  88.45  &  100  \\
0.6  &  100  &  100  &  99.5  &  100  \\
0.8  &  100  &  100  &  100  &  100  \\
1.0  &  100  &  100  &  100  &  100  \\ \hline
\end{tabular*}
\end{table}


%
\begin{table}
\tablewidth=10cm
\tabcolsep=0pt
\caption{Empirical power of test (in \%) based on 5000 simulations
using i.i.d. Brownian motions for $X_n(t)$, i.i.d. standard normals for
$\varepsilon_n$, $h(s,t)=cst$ and $N=500$}\label{t3}
\begin{tabular*}{10cm}{@{\extracolsep{\fill}}ld{2.2}d{2.2}d{2.2}d{2.2}d{2.2}@{}}
\hline
$c$&  \multicolumn{1}{l}{$p=1$} & \multicolumn{1}{l}{$p=2$} & \multicolumn{1}{l}{$p=3$} & \multicolumn{1}{l}{$p=4$} & \multicolumn{1}{l}{$p=5$}\\ \hline
& \multicolumn{5}{c}{$\alpha=0.01$} \\
0.0  &  1.05  &  1.25  &  0.8  &  1.35  &  1.85  \\
0.2  &  4.55  &  2.7  &  2.3  &  2.75  &  2.5  \\
0.4  &  19.4  &  11.55  &  8.2  &  7.4  &  5.9  \\
0.6  &  49.8  &  34.9  &  24.85  &  19.4  &  15.8  \\
0.8  &  79.75  &  64.5  &  50.95  &  41.95  &  36.55  \\
1.0  &  94.65  &  86.2  &  76.55  &  69.95  &  62.2  \\[5pt]
& \multicolumn{5}{c}{$\alpha=0.05$} \\
0.0  &  5.15  &  5.35  &  6.45  &  5.45  &  7.4  \\
0.2  &  13  &  8.75  &  8.15  &  9.4  &  9.8  \\
0.4  &  40.95  &  26.9  &  21.55  &  19.45  &  18.5  \\
0.6  &  71.25  &  58.2  &  46.2  &  39.3  &  35.3  \\
0.8  &  92.65  &  83.15  &  73.7  &  64.2  &  59  \\
1.0  &  98.65  &  94.85  &  89.75  &  86.2  &  80.3  \\[5pt]
& \multicolumn{5}{c}{$\alpha=0.10$} \\
0.0  &  10.15  &  10.15  &  11.95  &  11.95  &  13.45  \\
0.2  &  20.45  &  17.2  &  15.85  &  15.95  &  17.25  \\
0.4  &  53.55  &  38.95  &  32.15  &  28.7  &  29.05  \\
0.6  &  80.3  &  68.95  &  58.05  &  51.45  &  48.25  \\
0.8  &  95.95  &  90  &  83  &  75.7  &  69.9  \\
1.0  &  99.45  &  97.1  &  94.55  &  92.25  &  87.65  \\\hline
\end{tabular*}
\end{table}


%
\begin{table}
\tablewidth=11cm
\tabcolsep=0pt
\caption{Empirical power of test (in \%) based on 5000 simulations
using non-Gaussian $X_n(t)$ and non-normal $\varepsilon_n$,
$h(s,t)=cst$ and $N=500$}\label{t4}
\begin{tabular*}{11cm}{@{\extracolsep{\fill}}ld{2.2}d{2.2}d{2.2}d{2.2}@{}}
\hline
$c$&  \multicolumn{1}{l}{$p=1$} & \multicolumn{1}{l}{$p=2$} & \multicolumn{1}{l}{$p=3$} & \multicolumn{1}{l}{$p=4$} \\ \hline
& \multicolumn{4}{c}{$\alpha=0.01$} \\
0.0  &  1.65  &  1  &  1.4  &  2.15  \\
0.2  &  10.65  &  7.6  &  4.95  &  4.5  \\
0.4  &  35.45  &  32.3  &  25.7  &  22.2  \\
0.6  &  66.1  &  70.85  &  59.9  &  51.75  \\
0.8  &  87.85  &  90.65  &  87.15  &  81.25  \\
1.0  &  96.35  &  98.5  &  96.85  &  94.7  \\[5pt]
& \multicolumn{4}{c}{$\alpha=0.05$} \\
0.0  &  7.15  &  6.05  &  5.8  &  7.3  \\
0.2  &  23.25  &  18.9  &  15  &  15.55  \\
0.4  &  56.05  &  53.5  &  46.05  &  41.35  \\
0.6  &  81.6  &  84.05  &  76.5  &  69.25  \\
0.8  &  94.9  &  96.75  &  93.85  &  91.4  \\
1.0  &  98.45  &  99.6  &  99.05  &  98.25  \\[5pt]
& \multicolumn{4}{c}{$\alpha=0.10$} \\
0.0  &  13.45  &  11.3  &  11.1  &  14.6  \\
0.2  &  33.05  &  29.25  &  22.45  &  25.7  \\
0.4  &  66.5  &  63.5  &  56.8  &  52.9  \\
0.6  &  87.65  &  89.25  &  84  &  78  \\
0.8  &  97.25  &  98.4  &  95.8  &  94.45  \\
1.0  &  98.95  &  99.85  &  99.5  &  98.95  \\ \hline
\end{tabular*}
\end{table}


Although Theorem \ref{thmain} is proven under the assumption that
$X_n(t)$ is a Gaussian process, the result of Theorem \ref{thmain}
holds under relaxed conditions as discussed in Remark \ref{r1}. We
will now investigate the empirical size and power of our test when
$X_n(t)$ is not a Gaussian process. We generate the $\varepsilon_n$'s
according to a uniform distribution on $[-0.5,0.5]$. The predictors,
$X_n(t)$, are generated according to $X_n(t) = \xi_{n,1} + \xi_{n,2}t +
\xi_{n,3}(2t^2-1) + \xi_{n,4}(4t^3-3t)$, where $\{4\xi_{n,i}, 1\le i
\le4, 1\le n \}$ are i.i.d. random variables having a t-distribution
with 5 degrees of freedom. The polynomials in the definition of
$X_n(t)$ are the orthogonal Chebyshev polynomials. The resulting
empirical size and power are given in Table \ref{tviolations}. We see
from Table \ref{tviolations} that our testing procedure is robust
against non-Gaussian observations. Comparing Tables \ref{tempirical}
and \ref{tviolations}, we see that the value of the test statistics
tends to be larger if the $X_n$'s are not normally distributed for
small $N$. The overrejection fades as $N$ gets larger so in case of
non-Gaussian $X_n$'s, larger sample sizes are needed. This also
explains the somewhat better power of the procedure in the case of
non-Gaussian errors.

We also studied the choice of $p$ on the power of the test. The power
was studied under the alternative with the choice of $k(t)=1$,
$h(s,t)=c$ and $h(s,t)=cts$ using i.i.d. Brownian motions for $X_n$ and
i.i.d. standard normal errors for $\varepsilon_n$ with $N=500$ and
$p=1,\ldots,5$. We also repeated the simulations with $X_n$ chosen as
a $t_5$ process and $\varepsilon_n$ has a uniform distribution on
$[-0.5, 0.5]$. The results in Tables \ref{t1}--\ref{t4} illustrate that
choosing a larger $p$ might reduce the power of the test. We also
checked the power of the proposed test when $p=6,\ldots,12$ and the
result confirmed for these cases that the power is not necessarily a
monotone function of $p$ and using larger $p$'s might not provide
better testing method.

\section{Application to spectral data}\label{sapplication}

In this section, we apply our test to the data set collected by Tecator
and available at \url{http://lib.stat.cmu.edu/datasets/tecator}.
Tecator used 240 samples of finely chopped pure meat with different fat
contents. For each sample of meat, a 100 channel spectrum of
absorbances was recorded using a Tecator Infratec food and feed
analyzer. These absorbances can be thought of as a discrete
approximation to the continuous record, $X_n(t)$. Also, for each sample
of meat, the fat content, $Y_n$ was measured by analytic chemistry.

The absorbance curve measured from the $n${th} meat sample
is given by $X_n(t) = \log_{10}(I_0/I)$, where $t$ is the
wavelength of the light, $I_0$ is the intensity of the light before
passing through the meat sample, and $I$ is the intensity of the light
after it passes through the meat sample. The Tecator Infratec food and
feed analyzer measured absorbance at 100 different wavelengths between
850 and 1050 nanometers. This gives the values of $X_n(t)$ on a
discrete grid from which we can use cubic splines to interpolate the
values anywhere within the interval. A representative sample of 15 of
the 240 absorbance trajectories are pictured in Figure \ref{fspectra}.
Ferraty \textit{et al.} \cite{ferrvieupla2007} and Li and Yu \cite{liyu2008} contain an
analysis of the Tecator data as classification problem.

Yao and M\"uller \cite{muelleryao2010} proposed using a functional quadratic model
to predict the fat content, $Y_n$, of a meat sample based on its
absorbance spectrum, $X_n(t)$. We are interested in determining whether
the quadratic term in \eqref{eqmodel} is needed by testing its
significance for this data set. From the data, we calculate $U_{240}$.
The rejection probability is then $P(\chi^2(r) > U_{240})$.
The test statistic and hence the rejection probability are influenced
by the number of principal components that we choose to keep. If we
select $p$ according to the advice of Ramsay and Silverman \cite{ramsaysilverman2005},
we will keep only $p=1$ principal component because this explains more
than 85\% of the variation between absorbance curves in the sample.
Table \ref{tpvalues} gives rejection probabilities obtained using
$p=1$, $2$, and $3$ principal components, which strongly supports that
the quadratic regression provides a better model for the Tecator data.

\begin{table}[b]
\tablewidth=11cm
\tabcolsep=0pt
\caption[Rejection probabilities from absorbance data]{Rejection
probabilities (in \%) obtained by applying our testing procedure to the
Tecator data set with $p=1$, $2$, $3$, $4$, and $5$ principal components}\label{tpvalues}
\begin{tabular*}{11cm}{@{\extracolsep{\fill}}llllll@{}}
\hline
$p$ & $1$ & $2$ & $3$ & $4$ & $5$\\ \hline
Rejection probab. & $1.25$ & $13.15$ & $0.00$ & $0.00$ & $0.00$\\
\hline
\end{tabular*}
\end{table}

%
\section{\texorpdfstring{Proof of Theorem \protect\ref{thmain}}{Proof of Theorem 1.1}}\label{sproofs5}

We have from \eqref{eqconcise} and \eqref{eqhats} that
%
\begin{eqnarray}
\label{eqhats1} %
\pmatrix{ \hat{\mathbf A}
\cr
\hat{
\mathbf B}
\cr
\hat{\mu}} %
&=& \bigl(\hat{\mathbf Z}^T\hat{
\mathbf Z} \bigr)^{-1}\hat{\mathbf Z}^T \left(\hat{\mathbf
Z} %
\pmatrix{ \tilde{\mathbf A}
\cr
\tilde{\mathbf B}
\cr
\mu}
+ {\boldsymbol\varepsilon^{**}} \right)\nonumber
\\[-8pt]\\[-8pt]
&=& \pmatrix{ \tilde{\mathbf A}
\cr
\tilde{\mathbf B}
\cr
\mu} %
+
\bigl(\hat{\mathbf Z}^T\hat{\mathbf Z} \bigr)^{-1}\hat{
\mathbf Z}^T {\boldsymbol\varepsilon^{**}}.\nonumber %
\end{eqnarray}
We also note that, under the null hypothesis, $a_{i,j}=0$ for all $i$
and $j$ and therefore $\varepsilon_n^*$ and $\varepsilon_n^{**}$ of
\eqref{eqmodel1} and \eqref{eqmodel2} reduce to
\[
\varepsilon_n^*= \varepsilon_n + \sum
_{i=p+1}^{\infty} b_{i} \bigl\langle
X_n^c, v_i \bigr\rangle
\]
and
\[
\varepsilon_n^{**}=\varepsilon_n^* + \sum
_{i=1}^{p} b_{i} \bigl\langle
X_n^c, v_i-\hat{c}_i
\hat{v}_i \bigr\rangle+ \sum_{i=1}^p
b_i \langle\bar{X}-\mu_X,\hat{c}_i
\hat{v}_i \rangle.
\]
To obtain the limiting distribution of $\sqrt{N}\hat{\mathbf A}$, we
need to consider the vector\break $\sqrt{N} (\hat{\mathbf Z}^T\hat
{\mathbf Z} )^{-1}\hat{\mathbf Z}^T {\boldsymbol\varepsilon^{**}}$. We
will show in Lemmas \ref{lDD}--\ref{lC} that
%
\begin{equation}
\label{eqitsrighthere} \pmatrix{ \biggl(\dfrac{\hat{\mathbf Z}^T\hat{\mathbf Z}}{N} \biggr) - %
\pmatrix{ {\boldsymbol\zeta} {\mathbf G} {\boldsymbol\zeta} & {
\mathbf0}_{r\times p} & {\mathbf M}
\cr
{\mathbf0}_{p\times r} & {\boldsymbol
\Lambda}& {\mathbf 0}_{p\times1}
\cr
{\mathbf M}^T & {
\mathbf0}_{1\times p} & 1 } %
} = \mathrm{o}_{P} (1 ),
\end{equation}
where ${\boldsymbol\zeta}$ is an unobservable matrix of random signs,
${\boldsymbol\Lambda} = \diag(\lambda_1, \lambda_2, \dots,
\lambda_p)$, ${\mathbf M} = E ({\mathbf D}_n )$, and
${\mathbf G} = E ( {\mathbf D}_n {\mathbf D}_n^T  )$, where
\[
{\mathbf D}_n = \vech \bigl(\bigl\{ \bigl\langle v_i,
X_n^c \bigr\rangle \bigl\langle v_j ,
X_n^c \bigr\rangle , 1 \le i \le j \le p\bigr
\}^T \bigr). %
\]
We see from \eqref{eqitsrighthere} that the vector $\sqrt{N}
(\hat
{\mathbf Z}^T\hat{\mathbf Z} )^{-1}\hat{\mathbf Z}^T
{\boldsymbol
\varepsilon^{**}}$ has the same limiting distribution as
%
\begin{equation}
\label{eqhereagain} \frac{1}{\sqrt{N}} \sum_{n=1}^N
\varepsilon_n^{**} %
\pmatrix{ {\boldsymbol
\zeta} \bigl({\mathbf G}-{\mathbf M} {\mathbf M}^T
\bigr)^{-1} {\boldsymbol\zeta} & {\mathbf0}_{r\times p} & - {
\boldsymbol\zeta} \bigl({\mathbf G}-{\mathbf M} {\mathbf M}^T
\bigr)^{-1} {\boldsymbol\zeta} {\mathbf M}
\cr
{\mathbf0}_{p\times r}
& {\boldsymbol\Lambda}^{-1} & {\mathbf0}_{p\times1}
\cr
-{\mathbf
M}^T \bigl({\mathbf G}-{\mathbf M} {\mathbf M}^T
\bigr)^{-1} & {\mathbf0}_{1\times p} & 1+{\mathbf M}^T
\bigl({\mathbf G}-{\mathbf M} {\mathbf M}^T \bigr)^{-1} {
\mathbf M} } %
\pmatrix{ \hat{\mathbf D}_n
\cr
\hat{\mathbf
F}_n
\cr
1 } %
.
\end{equation}
Since we are only interested in $\sqrt{N} \hat{\mathbf A}$, we need
only consider the first $r=p(p+1)/2$ elements of the vector in \eqref
{eqhereagain}. In Lemma \ref{lestimatingG}, we show that these are
given by
\begin{eqnarray*}
&& \frac{1}{\sqrt{N}} \sum_{n=1}^N
\varepsilon_n^{**} %
\pmatrix{ {\boldsymbol\zeta}
\bigl({\mathbf G}-{\mathbf M} {\mathbf M}^T \bigr)^{-1} {
\boldsymbol\zeta} & {\mathbf0}_{r\times p} & - {\boldsymbol\zeta} \bigl({
\mathbf G}-{\mathbf M} {\mathbf M}^T \bigr)^{-1} {
\boldsymbol\zeta} {\mathbf M} } %
\pmatrix{ \hat{\mathbf
D}_n
\cr
\hat{\mathbf F}_n
\cr
1 }\\
%
&&\quad=
\frac{1}{\sqrt{N}} \sum_{n=1}^N
\varepsilon_n^{**} {\boldsymbol\zeta} \bigl({\mathbf G}-{
\mathbf M} {\mathbf M}^T \bigr)^{-1} {\boldsymbol\zeta} (
\hat{\mathbf D}_n - {\mathbf M} ).
\end{eqnarray*}
Then, in Lemma \ref{lthelimitingdistribution}, we prove that
\[
\frac{1}{\sqrt{N}} \sum_{n=1}^N
\varepsilon_n^{**} \bigl({\mathbf G}-{\mathbf M} {\mathbf
M}^T \bigr)^{-1} {\boldsymbol\zeta} (\hat {\mathbf
D}_n - {\mathbf M} ) \inD N \bigl(0, \tau^2 \bigl({
\mathbf G}-{\mathbf M} {\mathbf M}^T \bigr)^{-1} \bigr),
\]
where $\tau^2 = \var (\varepsilon_1^* )$.
Finally, in Lemmas \ref{lestimatorsconverging} and \ref{ltau}, we
show that $\hat{\tau}^2 - \tau^2 = \mathrm{o}_{P} (1 )$. As a
consequence of \eqref
{eqitsrighthere}, we see that $ (\hat{\mathbf G}- \hat{\mathbf
M}\hat{\mathbf M}^T ) - {\boldsymbol\zeta} ({\mathbf G}-
{\mathbf M}{\mathbf M}^T ){\boldsymbol\zeta}=\mathrm{o}_{P}
(1 )$. Since
${\boldsymbol\zeta}$ is a diagonal matrix of signs, ${\boldsymbol
\zeta
}{\boldsymbol\zeta} = I$, completing the proof of Theorem \ref{thmain}.

\section{\texorpdfstring{Proof of Theorem \protect\ref{thconst}}{Proof of Theorem 1.2}}\label{sconsistency}

We provide only an outline of the proof since it follows the arguments
used in the proof of Theorem \ref{thmain}. However, the arguments are
simple since instead of obtaining an asymptotic limit distribution we
only establish the weak law
%
\begin{equation}
\label{eqprobab} \hat{\mathbf A}^T\bigl(\hat{{\mathbf G}}-\hat{{
\mathbf M}}\hat{{\mathbf M}}^T\bigr)\hat{\mathbf A} \inP{\mathbf
A}^T\bigl({{\mathbf G}}-{{\mathbf M}} {{\mathbf M}}^T
\bigr){\mathbf A},
\end{equation}
where ${\mathbf A} = \vech ( \{ a_{i,j} (2-1\{i=j\} ),
1\le i \le j\le p \}^T )$ is like the vector $\tilde{\mathbf A}$
except without the random signs.

First, we note that according to Lemma \ref{lgoodapproximation}, the
estimation of $v_1,\ldots, v_p$ by $\hat{v}_1,\ldots, \hat{v}_p$
causes only the \hyperref[sintro5]{Introduction} of the random signs $\hat{c}_1, \ldots
,\hat{c}_p$. As in the proof of Theorem \ref{thmain}, one can verify that
\[
\hat{\mathbf A}-{\boldsymbol\zeta} {\mathbf A} \inP{\mathbf0}.
\]
Lemmas \ref{lDD} and \ref{lD} hold under $H_0$ as well as under
$H_A$. This gives
\[
\hat{{\mathbf G}}-{\boldsymbol\zeta} {{\mathbf G}} {\boldsymbol
\zeta}=\mathrm{o}_P(1)
\]
and
\[
\hat{{\mathbf M}}\hat{{\mathbf M}}^T-{\boldsymbol\zeta} {{\mathbf
M}} {{\mathbf M}}^T{\boldsymbol\zeta} =\mathrm{o}_P(1),
\]
completing the proof of \eqref{eqprobab}.

\section{Technical lemmas}\label{stech}

Throughout the proofs in this section, we will use $\|\cdot\|_1$ to be
the 1-norm and $\|\cdot\|_2$ to be 2-norm on the unit interval, square,
cube, or hypercube. The null hypothesis, $H_0$, is assumed throughout
this section. We will make frequent use of the following lemma, which
is established in Dauxois \textit{et al.} \cite{dauxois1982} and Bosq \cite{bosq2000}.

\begin{lemma}\label{lgoodapproximation}
If Assumptions \ref{agaussian}, \ref{aXmoments}, and \ref{aunique}
hold, then
\[
\bigl\|\hat{c}_i\hat{v}_i(t) - v_i(t)
\bigr\|_2 = \mathrm{O}_{P}\bigl(N^{-1/2}\bigr)
\]
for each $1\le i \le p$.
\end{lemma}
%
\begin{lemma}\label{lDD}
If Assumptions \ref{agaussian}, \ref{aXmoments}, and \ref{aunique}
hold, then there is a non-random matrix ${\mathbf G}$ such that
\[
(\hat{\mathbf G} - {\boldsymbol\zeta} {\mathbf G} {\boldsymbol \zeta} ) =
\mathrm{o}_{P} (1 ),
\]
where $\hat{\mathbf G} = N^{-1} \sum_{n=1}^{N} \hat{\mathbf D}_n
\hat
{\mathbf D}_n^T$ and ${\boldsymbol\zeta} = \diag (\vech(\{\hat
{c}_i \hat{c}_j, 1\le i \le j \le p\}^T) )$.
\end{lemma}
\begin{pf}
By the Karhunen--Lo\'eve expansion, we have
%
\begin{equation}
\label{eqkarlo} X_n^c(t)=\sum
_{\ell=1}^\infty\lambda^{1/2}_\ell
\xi_\ell^{(n)} v_\ell(t).
\end{equation}
Therefore an element of ${\mathbf D}_n {\mathbf D}_n^T$ is of the form
$\sqrt{\lambda_i \lambda_j \lambda_k \lambda_{\ell}} \xi_i^{(n)}
\xi_j^{(n)} \xi_k^{(n)} \xi_{\ell}^{(n)}$. Hence, using the strong law of
large numbers we conclude
\[
\frac{1}{N}\sum_{n=1}^{N} {\mathbf
D}_n {\mathbf D}_n^T \as{\mathbf G},
\]
where ${\mathbf G} = E ( {\mathbf D}_n {\mathbf D}_n^T  )$.
Thus, it suffices to show that
%
\begin{equation}
\label{eqmatrixwise} \frac{1}{N} \sum_{n=1}^{N}
\bigl({\boldsymbol\zeta} \hat{\mathbf D}_n \hat{\mathbf
D}_n^T {\boldsymbol\zeta} - {\mathbf D}_n {
\mathbf D}_n^T \bigr) = \mathrm{o}_{P} (1 ).
\end{equation}
Expressing \eqref{eqmatrixwise} elementwise, we obtain
%
\begin{eqnarray}
\label{eqgammabreakup} &&\frac{1}{N} \sum_{n=1}^{N}
\bigl(\langle X_n- \bar{X}, \hat {c}_i
\hat{v}_i \rangle \langle X_n- \bar{X},
\hat{c}_j \hat{v}_j \rangle \langle X_n-
\bar{X}, \hat{c}_k \hat{v}_k \rangle \langle
X_n- \bar{X}, \hat {c}_{\ell} \hat {v}_{\ell}
\rangle
\nonumber
\\[-8pt]
\\[-8pt]
&&\hspace*{30pt}{}- \bigl\langle X_n^c, v_i \bigr\rangle
\bigl\langle X_n^c, v_j \bigr\rangle \bigl
\langle X_n^c, v_k \bigr\rangle \bigl\langle
X_n^c, v_{\ell} \bigr\rangle
\bigr)=\mathrm{o}_{P} (1 ).
\nonumber
\end{eqnarray}
In order to prove \eqref{eqgammabreakup}, it is enough to show that
%
\begin{eqnarray}
\label{eqfavorite} %
&&\frac{1}{N} \sum
_{n=1}^{N} \bigl(\bigl\langle X_n^c,
\hat{c}_i \hat {v}_i \bigr\rangle \bigl\langle
X_n^c, \hat{c}_j \hat{v}_j
\bigr\rangle \bigl\langle X_n^c, \hat{c}_k
\hat {v}_k \bigr\rangle \bigl\langle X_n^c,
\hat{c}_{\ell} \hat{v}_{\ell} \bigr\rangle
\nonumber
\\[-8pt]
\\[-8pt]
&&\hspace*{30pt}{}- \bigl\langle X_n^c, v_i \bigr\rangle
\bigl\langle X_n^c, v_j \bigr\rangle \bigl
\langle X_n^c, v_k \bigr\rangle \bigl\langle
X_n^c, v_{\ell} \bigr\rangle
\bigr)=\mathrm{o}_{P} (1 )
\nonumber
\end{eqnarray}
and
%
\begin{eqnarray}
\label{eqnextbestthing} %
&&\frac{1}{N} \sum
_{n=1}^{N} \bigl(\langle X_n- \bar{X},
\hat {c}_i \hat{v}_i \rangle \langle X_n-
\bar{X}, \hat{c}_j \hat{v}_j \rangle \langle
X_n- \bar{X}, \hat{c}_k \hat{v}_k \rangle
\langle X_n- \bar{X}, \hat {c}_{\ell} \hat {v}_{\ell}
\rangle
\nonumber
\\[-8pt]
\\[-8pt]
&&\hspace*{30pt}{}- \bigl\langle X_n^c, \hat{c}_i
\hat{v}_i \bigr\rangle \bigl\langle X_n^c,
\hat{c}_j \hat{v}_j \bigr\rangle \bigl\langle
X_n^c, \hat{c}_k\hat{v}_k
\bigr\rangle \bigl\langle X_n^c, \hat {c}_{\ell
}
\hat{v}_{\ell} \bigr\rangle \bigr)=\mathrm{o}_{P} (1 ).
\nonumber
\end{eqnarray}

We only establish \eqref{eqfavorite}, since the proof of \eqref
{eqnextbestthing} is essentially the same. Using H\"older's
inequality, we obtain
\begin{eqnarray*}
&& \Biggl|\Int\Int\Int\Int \Biggl(\frac{1}{N} \sum_{n=1}^{N}
X_n^c(s)X_n^c(t)X_n^c(u)X_n^c(w)
\Biggr)
\\
&&\hspace*{62pt}{}\times \bigl(\hat{c}_i \hat{v}_i(s)
\hat{c}_j \hat{v}_j(t) \hat{c}_k
\hat{v}_k(u) \hat{c}_{\ell} \hat{v}_{\ell}(w) -
v_i(s) v_j(t) v_k(u) v_{\ell}(w)
\bigr)\,\mathrm{d}s\,\mathrm{d}t\,\mathrm{d}u\,\mathrm{d}w \Biggr|
\\
&&\quad\le\Biggl\Vert  \frac{1}{N} \sum
_{n=1}^{N} X_n^c(s)X_n^c(t)X_n^c(u)X_n^c(w)
 \Biggr\Vert_2
\\
&&\qquad{}\times \bigl\Vert \hat{c}_i \hat{v}_i(s)
\hat{c}_j \hat {v}_j(t) \hat{c}_k
\hat{v}_k(u) \hat{c}_{\ell} \hat{v}_{\ell}(w) -
v_i(s) v_j(t) v_k(u) v_{\ell}(w)
\bigr\Vert_2.
\end{eqnarray*}
By the law of large numbers in Hilbert spaces (cf. Bosq \cite{bosq2000}),
we have that
\[
 \Biggl\Vert \frac{1}{N} \sum_{n=1}^{N}
X_n^c(s)X_n^c(t)X_n^c(u)X_n^c(w)
\Biggr\Vert_2=\mathrm{O}_{P}(1),
\]
so it remains only to show that
\[
\bigl\Vert \hat{c}_i \hat{v}_i(s)
\hat{c}_j \hat{v}_j(t) \hat {c}_k \hat
{v}_k(u) \hat{c}_{\ell} \hat{v}_{\ell}(w) -
v_i(s) v_j(t) v_k(u) v_{\ell
}(w)
 \bigr\Vert_2 = \mathrm{o}_{P} (1 ).
\]
Using Minkowski's inequality, Fubini's theorem, the fact that $\|\hat
{v}_{i}\|_2 = \|{v}_{i}\|_2 = 1$, and then Lemma \ref
{lgoodapproximation}, we obtain
\begin{eqnarray*}
&&\bigl\Vert \hat{c}_i \hat{v}_i(s)
\hat{c}_j \hat{v}_j(t) \hat{c}_k
\hat{v}_k(u) \hat{c}_{\ell} \hat{v}_{\ell}(w) -
v_i(s) v_j(t) v_k(u) v_{\ell}(w)
 \bigr\Vert_2
\\
&&\quad\le\bigl\Vert \bigl(\hat{c}_i \hat{v}_i(s)
- v_i(s) \bigr) \hat{c}_j \hat{v}_j(t)
\hat{c}_k \hat{v}_k(u) \hat{c}_{\ell} \hat
{v}_{\ell}(w)  \bigr\Vert_2
\\
&&\qquad{}+ \bigl\Vert {v}_i(s) \hat{c}_j
\hat{v}_j(t) \hat{c}_k \hat{v}_k(u) \bigl(
\hat{c}_{\ell} \hat{v}_{\ell
}(w)-v_{\ell
}(w) \bigr)
\bigr\Vert_2
\\
&&\qquad{}+ \bigl\Vert {v}_i(s) \hat{c}_j
\hat{v}_j(t) \bigl(\hat{c}_k \hat{v}_k(u)-
v_k(u) \bigr) {v}_{\ell}(w)  \bigr\Vert_2
\\
&&\qquad{}+ \bigl\Vert {v}_i(s) \bigl(\hat{c}_j
\hat {v}_j(t)-v_j(t) \bigr) {v}_k(u)
{v}_{\ell}(w)  \bigr\Vert_2
\\
&&\quad=\bigl\Vert \hat{c}_i \hat{v}_i(s)-v_i(s)
 \bigr\Vert_2 + \bigl\Vert
\hat{c}_j \hat{v}_j(t)-v_j(t) \bigr\Vert_2 + \bigl\Vert \hat{c}_k
\hat{v}_k(u)-v_k(u) \bigr\Vert_2 +\bigl\Vert \hat {c}_{\ell}
\hat{v}_{\ell}(w)-v_{\ell}(w) \bigr\Vert_2
\\
&&\quad = \mathrm{O}_{P}\bigl(N^{-1/2}\bigr).
\end{eqnarray*}
Hence, \eqref{eqfavorite} is proven which also completes the proof of
Lemma \ref{lDD}.
\end{pf}

\begin{lemma}\label{lFD}
If Assumptions \ref{agaussian}, \ref{aXmoments}, and \ref{aunique}
hold, then
\[
\frac{1}{N} \sum_{n=1}^{N} \hat{
\mathbf F}_n \hat{\mathbf D}_n^T
=\mathrm{o}_{P} (1 ).
\]
\end{lemma}
\begin{pf}
We see from \eqref{eqkarlo} that an element of ${\mathbf F}_n {\mathbf
D}_n^T$ can be written in the form $\sqrt{\lambda_i \lambda_j
\lambda_k} \xi_i^{(n)} \xi_j^{(n)}\* \xi_k^{(n)}$, where ${\mathbf F}_n =
(\langle X_n^c, v_1 \rangle,
\langle X_n^c, v_2 \rangle,
\dots,
\langle X_n^c, v_p \rangle)
^T$. We observe that
$E\xi_i^{(n)} \xi_j^{(n)} \xi_k^{(n)}=0$, so using the central limit
theorem, we have
\[
\frac{1}{N} \sum_{n=1}^{N} {\mathbf
F}_n {\mathbf D}_n^T = \mathrm{O}_{P}
\bigl(N^{-1/2}\bigr).
\]
Repeating the arguments in the proof \eqref{eqgammabreakup}, one can
verify that
%
\begin{eqnarray}
\label{eqbus}
&&\frac{1}{N} \sum_{n=1}^{N}
\bigl(\langle X_n-\bar{X}, \hat{c}_i
\hat{v}_i \rangle \langle X_n-\bar{X},
\hat{c}_j \hat{v}_j \rangle \langle X_n-
\bar{X}, \hat {c}_k \hat{v}_k \rangle
\nonumber
\\[-8pt]
\\[-8pt]
&&\hspace*{30pt}{}- \bigl\langle X_n^c, v_i \bigr\rangle
\bigl\langle X_n^c, v_j \bigr\rangle \bigl
\langle X_n^c, v_k \bigr\rangle \bigr) =
\mathrm{o}_{P} (1 ).
\nonumber
\end{eqnarray}
Since random signs do not affect convergence to zero, the proof is complete.
\end{pf}

\begin{lemma}\label{lFF}
If Assumptions \ref{agaussian}, \ref{aXmoments}, and \ref{aunique}
hold, then
\[
\frac{1}{N}\sum_{n=1}^{N} \hat{
\mathbf F}_n \hat{\mathbf F}_n^T - {
\boldsymbol\Lambda} = \mathrm{o}_{P} (1 ),
\]
where ${\boldsymbol\Lambda} = {\diag}(\lambda_1, \lambda_2, \dots,
\lambda_p)$.
\end{lemma}
\begin{pf}
By \eqref{eqkarlo}, an element of ${\mathbf F}_n {\mathbf F}_n^T$ is
of the form $\sqrt{\lambda_i \lambda_j } \xi_i^{(n)} \xi_j^{(n)}$.
Since $E\xi_i^{(n)} \xi_j^{(n)}=1\{i=j\}$, according to the law of
large numbers, we have
\[
\frac{1}{N} \sum_{n=1}^{N} {\mathbf
F}_n {\mathbf F}_n^T - {\boldsymbol
\Lambda}= \mathrm{o}_{P} (1 ).
\]
Thus, it suffices to demonstrate that
%
\begin{equation}
\label{eqbumpy} \frac{1}{N} \sum_{n=1}^{N}
\bigl(\langle X_n-\bar{X}, \hat{v}_i \rangle \langle
X_n-\bar{X}, \hat{v}_j \rangle- \bigl\langle
X_n^c, v_i \bigr\rangle \bigl\langle
X_n^c, v_j \bigr\rangle \bigr) =
\mathrm{o}_{P} (1 ).
\end{equation}
Since random signs do not affect convergence to zero, multiplying $\hat
{v}_i$ by $\hat{c}_i$ and $\hat{v}_j$ by $\hat{c}_j$ will not affect
convergence when $i\neq j$. If $i=j$, then $\hat{c}_i\hat{c}_j = \hat
{c}_i^2 = 1$. Therefore, it suffices to show that
%
\begin{equation}
\label{eqbumpyride} \frac{1}{N} \sum_{n=1}^{N}
\bigl(\langle X_n-\bar{X}, \hat{c}_i
\hat{v}_i \rangle \langle X_n-\bar{X},
\hat{c}_j \hat{v}_j \rangle- \bigl\langle
X_n^c, v_i \bigr\rangle \bigl\langle
X_n^c, v_j \bigr\rangle \bigr) =
\mathrm{o}_{P} (1 ).
\end{equation}
One can show \eqref{eqbumpyride} in exactly the same way we
established \eqref{eqgammabreakup} in the proof of Lemma \ref{lDD}.
This completes the proof.
\end{pf}

\begin{lemma}\label{lF}
If Assumptions \ref{agaussian}, \ref{aXmoments}, and \ref{aunique}
hold, then
\[
\frac{1}{N}\sum_{n=1}^{N} \hat{
\mathbf F}_n = \mathrm{o}_{P} (1 ).
\]
\end{lemma}
\begin{pf}
Using \eqref{eqkarlo}, an element of ${\mathbf F}_n$ has the form
$\sqrt{\lambda_i}\xi_i^{(n)}$, so the law of large numbers implies that
\[
\frac{1}{N}\sum_{n=1}^{N} {\mathbf
F}_n = \mathrm{o}_{P} (1 ).
\]
The proof will be completed by establishing that
%
\begin{equation}
\label{eqtheforgotten} \frac{1}{N}\sum_{n=1}^{N}
({\mathbf F}_n -\hat{\mathbf F}_n ) =\mathrm{o}_{P}
(1 ).
\end{equation}
We express \eqref{eqtheforgotten} componentwise and obtain
%
\begin{equation}
\label{eqtheforgotten2} \frac{1}{N}\sum_{n=1}^{N}
\bigl(\bigl\langle X_n^c, v_i\bigr\rangle-
\langle X_n - \bar{X}, \hat{v}_i \rangle \bigr)
=\mathrm{o}_{P} (1 ).
\end{equation}
Since random signs do not affect convergence to zero, it suffices to
show that
%
\begin{equation}
\label{eqtheforgotten3} \frac{1}{N}\sum_{n=1}^{N}
\bigl(\bigl\langle X_n^c, v_i\bigr\rangle-
\langle X_n - \bar{X}, \hat{c}_i \hat{v}_i
\rangle \bigr) =\mathrm{o}_{P} (1 ).
\end{equation}
We will establish \eqref{eqtheforgotten3} in two steps. We will show that
%
\begin{equation}
\label{eqpeanut} \frac{1}{N}\sum_{n=1}^{N}
\bigl(\bigl\langle X_n^c, v_i\bigr\rangle-
\bigl\langle X_n^c, \hat {c}_i
\hat{v}_i \bigr\rangle \bigr) =\mathrm{o}_{P} (1 ).
\end{equation}
Then, we will establish that
%
\begin{equation}
\label{eqbutter} \frac{1}{N}\sum_{n=1}^{N}
\bigl(\bigl\langle X_n^c, \hat{c}_i \hat
{v}_i\bigr\rangle - \langle X_n - \bar{X},
\hat{c}_i \hat{v}_i \rangle \bigr) =\mathrm{o}_{P} (1
).
\end{equation}

Using the central limit theorem in Hilbert spaces with Lemma \ref
{lgoodapproximation}, we conclude
\begin{eqnarray*}
\Biggl\llvert \frac{1}{N}\sum_{n=1}^{N}
\bigl(\bigl\langle X_n^c, v_i\bigr\rangle-
\bigl\langle X_n^c, \hat{c}_i
\hat{v}_i \bigr\rangle \bigr)\Biggr\rrvert &\le&\Biggl
\Vert \frac{1}{N}\sum_{n=1}^{N}
X_n^c(t) \bigl(v_i(t)-\hat
{c}_i \hat{v}_i(t) \bigr) \Biggr
\Vert_1
\\
&\le&\Biggl\Vert  \frac{1}{N}\sum
_{n=1}^{N} X_n^c(t) \Biggr
\Vert_2 \bigl\Vert v_i(t)-
\hat{c}_i \hat{v}_i (t)  \bigr\Vert_2
\\
&=&\mathrm{O}_{P}\bigl(N^{-1}\bigr),
\end{eqnarray*}
and by the same arguments we have
\begin{eqnarray*}
\Biggl\llvert \frac{1}{N}\sum_{n=1}^{N}
\bigl(\bigl\langle X_n^c, \hat{c}_i \hat
{v}_i\bigr\rangle - \langle X_n - \bar{X},
\hat{c}_i \hat{v}_i \rangle \bigr)\Biggr\rrvert
&=&\bigl\llvert \langle\mu_X - \bar{X},
\hat{c}_i \hat{v}_i\rangle \bigr\rrvert
\\
&\le&\bigl\Vert \bigl(\mu_X(t) - \bar{X}(t) \bigr)
\hat{c}_i \hat {v}_i(t)  \bigr\Vert_1
\\
&\le&\bigl\Vert \mu_X(t) - \bar{X}(t) \bigr\Vert_2
\\
&=&\mathrm{o}_{P} (1 ).
\end{eqnarray*}
\upqed\end{pf}

\begin{lemma}\label{lD}
If Assumptions \ref{agaussian}, \ref{aXmoments}, and \ref{aunique}
hold, then
\[
\hat{\mathbf M} - {\mathbf M} = \mathrm{o}_{P} (1 ).
\]
where $\hat{\mathbf M} = N^{-1}\sum_{n=1}^{N} \hat{\mathbf D}_n$ and
${\mathbf M} = E ({\mathbf D}_n )$.
\end{lemma}
\begin{pf}
An arbitrary element of $\hat{\mathbf D}_n$ is of the form
\[
\frac{1}{N} \sum_{n=1}^{N} \langle
X_n-\bar{X}, \hat{v}_i \rangle \langle
X_n-\bar {X}, \hat{v}_j \rangle.
\]
Since this is exactly the same as the form of an arbitrary element of
$\hat{\mathbf F}_n \hat{\mathbf F}_n^T$, Lemma \ref{lD} follows from
the proof of Lemma \ref{lFF}. Note in particular that when $i\neq j$,
the sum converges to zero and is unaffected by signs, and when $i=j$,
the signs cancel each other out. For this reason, ${\boldsymbol\zeta}
{\mathbf M} = {\mathbf M}$, rendering it unnecessary to multiply
${\mathbf M}$ by ${\boldsymbol\zeta}$ in the statement of the lemma.
\end{pf}

\begin{lemma}\label{lC}
If Assumptions \ref{agaussian}, \ref{aXmoments}, and \ref{aunique}
hold, then
\[
\pmatrix{ \biggl(\dfrac{\hat{\mathbf Z}^T\hat{\mathbf Z}}{N} \biggr) - %
\pmatrix{ {\boldsymbol
\zeta} {\mathbf G} {\boldsymbol\zeta} & {\mathbf0}_{r\times p} & {\mathbf M}
\cr
{\mathbf0}_{p\times r} & {\boldsymbol\Lambda}& {\mathbf 0}_{p\times1}
\cr
{
\mathbf M}^T & {\mathbf0}_{1\times p} & 1 } %
} =
\mathrm{o}_{P} (1 ).
\]
\end{lemma}
\begin{pf}
This follows immediately from Lemmas \ref{lDD}--\ref{lD}.
\end{pf}

We will now use Lemma \ref{lC} to separate our estimate, $\hat
{\mathbf
A}$, of $\tilde{\mathbf A}$ from the estimates of the other parameters
in \eqref{eqhats}.
%
\begin{lemma}\label{lestimatingG}
If Assumptions \ref{agaussian}--\ref{aunique} hold, then
\[
{\boldsymbol\zeta} \sqrt{N} \hat{\mathbf A} - N^{-1/2} \sum
_{n=1}^N \varepsilon_n^{**}
\bigl({\mathbf G}-{\mathbf M} {\mathbf M}^T \bigr)^{-1} {
\boldsymbol\zeta} (\hat{\mathbf D}_n - {\mathbf M} ) =
\mathrm{o}_{P} (1 ).
\]
\end{lemma}
\begin{pf}
Let
\[
{\mathbf C} = %
\pmatrix{ {\boldsymbol\zeta} {\mathbf G} {\boldsymbol
\zeta} & {\mathbf0}_{r\times p} & {\mathbf M}
\cr
{\mathbf0}_{p\times r} & {
\boldsymbol\Lambda}& {\mathbf 0}_{p\times1}
\cr
{\mathbf M}^T & {
\mathbf0}_{1\times p} & 1 }.
\]
Using the fact that ${\boldsymbol\zeta}{\mathbf M} = {\mathbf M}$, one
can verify via matrix multiplication that
\[
{\mathbf C}^{-1} = %
\pmatrix{ {\boldsymbol\zeta} \bigl({
\mathbf G}-{\mathbf M} {\mathbf M}^T \bigr)^{-1} {
\boldsymbol\zeta} & {\mathbf0}_{r\times p} & - {\boldsymbol\zeta} \bigl({
\mathbf G}-{\mathbf M} {\mathbf M}^T \bigr)^{-1} {
\boldsymbol\zeta} {\mathbf M}
\cr
{\mathbf0}_{p\times r} & {\boldsymbol
\Lambda}^{-1} & {\mathbf0}_{p\times1}
\cr
-{\mathbf M}^T
\bigl({\mathbf G}-{\mathbf M} {\mathbf M}^T \bigr)^{-1} & {
\mathbf0}_{1\times p} & 1+{\mathbf M}^T \bigl({\mathbf G}-{\mathbf
M} {\mathbf M}^T \bigr)^{-1} {\mathbf M} }.
\]
Since $N^{-1/2}\hat{\mathbf Z}^T{\boldsymbol\varepsilon^{**}}$ %
is bounded in probability, by \eqref{eqhats1} and Lemma \ref{lC} we have
%
\begin{equation}
\label{eqthehatsareback} \sqrt{N} %
\pmatrix{ \hat{\mathbf A}
\cr
\hat{
\mathbf B} - \tilde{\mathbf B}
\cr
\hat{\mu} - \mu } %
- {\mathbf
C}^{-1} N^{-1/2}\hat{\mathbf Z}^T{\boldsymbol
\varepsilon^{**}} = \mathrm{o}_{P} (1 ).
\end{equation}
We observe that ${\mathbf C}^{-1} N^{-1/2}\hat{\mathbf
Z}^T{\boldsymbol
\varepsilon^{**}}$ can be expressed as
%
\begin{equation}
\label{eqrowsonly} N^{-1/2} \sum_{n=1}^N
\varepsilon_n^{**} %
\pmatrix{ {\boldsymbol
\zeta} \bigl({\mathbf G}-{\mathbf M} {\mathbf M}^T
\bigr)^{-1} {\boldsymbol\zeta} & {\mathbf0}_{r\times p} & - {
\boldsymbol\zeta} \bigl({\mathbf G}-{\mathbf M} {\mathbf M}^T
\bigr)^{-1} {\boldsymbol\zeta} {\mathbf M}
\cr
{\mathbf0}_{p\times r}
& {\boldsymbol\Lambda}^{-1} & {\mathbf0}_{p\times1}
\cr
-{\mathbf
M}^T \bigl({\mathbf G}-{\mathbf M} {\mathbf M}^T
\bigr)^{-1} & {\mathbf0}_{1\times p} & 1+{\mathbf M}^T
\bigl({\mathbf G}-{\mathbf M} {\mathbf M}^T \bigr)^{-1} {
\mathbf M} } %
\pmatrix{ \hat{\mathbf D}_n
\cr
\hat{\mathbf
F}_n
\cr
1}.
\end{equation}
Notice that the first $r=p(p+1)/2$ elements of the vector in \eqref
{eqrowsonly} are given by
\begin{eqnarray*}
&& N^{-1/2} \sum_{n=1}^N
\varepsilon_n^{**} %
\pmatrix{ {\boldsymbol\zeta}
\bigl({\mathbf G}-{\mathbf M} {\mathbf M}^T \bigr)^{-1} {
\boldsymbol\zeta} & {\mathbf0}_{r\times p} & - {\boldsymbol\zeta} \bigl({
\mathbf G}-{\mathbf M} {\mathbf M}^T \bigr)^{-1} {
\boldsymbol\zeta} {\mathbf M} } %
\pmatrix{ \hat{\mathbf
D}_n
\cr
\hat{\mathbf F}_n
\cr
1} %
\\
&&\quad= N^{-1/2} \sum_{n=1}^N
\varepsilon_n^{**} \bigl( {\boldsymbol\zeta} \bigl({\mathbf
G}-{\mathbf M} {\mathbf M}^T \bigr)^{-1} {\boldsymbol\zeta}
\hat{\mathbf D}_n - {\boldsymbol\zeta} \bigl({\mathbf G}-{\mathbf M}
{\mathbf M}^T \bigr)^{-1} {\boldsymbol\zeta} {\mathbf M}
\bigr)
\\
&&\quad= N^{-1/2} \sum_{n=1}^N
\varepsilon_n^{**} {\boldsymbol \zeta} \bigl({\mathbf G}-{
\mathbf M} {\mathbf M}^T \bigr)^{-1} {\boldsymbol\zeta} (
\hat{\mathbf D}_n - {\mathbf M} ).
\\
\end{eqnarray*}
Therefore,
%
\begin{equation}
\label{eqbabysteps} \sqrt{N} \hat{\mathbf A} - N^{-1/2} \sum
_{n=1}^N \varepsilon_n^{**} {
\boldsymbol\zeta} \bigl({\mathbf G}-{\mathbf M} {\mathbf M}^T
\bigr)^{-1} {\boldsymbol\zeta} (\hat{\mathbf D}_n - {
\mathbf M} ) = \mathrm{o}_{P} (1 ).
\end{equation}
The result is now obtained by multiplying \eqref{eqbabysteps} on the
left by ${\boldsymbol\zeta}$.
\end{pf}

\begin{lemma}\label{lthelimitingdistribution}
If Assumptions \ref{agaussian}--\ref{aunique} hold, then
\[
N^{-1/2} \sum_{n=1}^N
\varepsilon_n^{**} \bigl({\mathbf G}-{\mathbf M} {\mathbf
M}^T \bigr)^{-1} {\boldsymbol\zeta} (\hat{\mathbf
D}_n - {\mathbf M} ) \inD N \bigl(0, \tau^2 \bigl({
\mathbf G}-{\mathbf M} {\mathbf M}^T \bigr)^{-1} \bigr),
\]
where
\[
\tau^2 = \sigma^2 + \sum_{i=p+1}^{\infty}b_i^2
\lambda_i
\]
and $\sigma^2=\var\varepsilon_n$.
\end{lemma}
\begin{pf}
We prove this lemma in three steps. First, we establish that
%
\begin{equation}
\label{eqtofu} N^{-1/2} \sum_{n=1}^N
\varepsilon_n^{**} \bigl( ({\boldsymbol \zeta }\hat{
\mathbf D}_n - {\mathbf M} ) - ( {\mathbf D}_n - {
\mathbf M} ) \bigr) = \mathrm{o}_{P} (1 ).
\end{equation}
In the second step, we prove that
%
\begin{equation}
\label{eqicecream} N^{-1/2} \sum_{n=1}^N
( {\mathbf D}_n - {\mathbf M} ) \Biggl(\varepsilon_n^{**}
- \varepsilon_n^{*} - \sum_{i=1}^p
b_i \langle \bar{X} - \mu_X , \hat{c}_i
\hat{v}_i \rangle \Biggr) = \mathrm{o}_{P} (1 )
\end{equation}
and
%
\begin{equation}
\label{eqorangejuice} N^{-1/2} \sum_{n=1}^N
( {\mathbf D}_n - {\mathbf M} )\langle \bar {X} - \mu_X
, \hat{c}_i \hat{v}_i \rangle= \mathrm{o}_{P} (1 ).
\end{equation}
Combining \eqref{eqtofu}, \eqref{eqicecream}, and \eqref
{eqorangejuice}, we obtain immediately that
\[
N^{-1/2} \sum_{n=1}^N \bigl({
\mathbf G}-{\mathbf M} {\mathbf M}^T \bigr)^{-1} \bigl(
\varepsilon_n^{**} ({\boldsymbol\zeta} \hat {\mathbf
D}_n - {\mathbf M} ) - \varepsilon_n^{*} ({
\mathbf D}_n - {\mathbf M} ) \bigr) = \mathrm{o}_{P} (1 ).
\]
Therefore, the lemma will be established by the third step:
%
\begin{equation}
\label{eqredpen} N^{-1/2} \sum_{n=1}^N
\bigl({\mathbf G}-{\mathbf M} {\mathbf M}^T \bigr)^{-1}
\varepsilon_n^{*} ({\mathbf D}_n - {\mathbf
M} ) \inD N \bigl(0, \tau^2 \bigl({\mathbf G}-{\mathbf M} {\mathbf
M}^T \bigr)^{-1} \bigr).
\end{equation}

We will now proceed to prove \eqref{eqtofu}. The left-hand side of \eqref
{eqtofu} can be expressed elementwise as
%
\begin{equation}
\label{eqrunner} N^{-1/2} \sum_{n=1}^N
\varepsilon_n^{**} \bigl(\langle X_n - \bar
{X} , \hat {c}_i \hat{v}_i \rangle \langle
X_n - \bar{X} , \hat{c}_j \hat{v}_j
\rangle- \bigl\langle X_n^c , v_i \bigr
\rangle \bigl\langle X_n^c , v_j \bigr
\rangle \bigr) = \mathrm{o}_{P} (1 ),
\\
\end{equation}
so it is sufficient to show that
%
\begin{equation}
\label{eqtablespoon} N^{-1/2} \sum_{n=1}^N
\varepsilon_n^{**} \bigl(\bigl\langle X_n^c
, \hat{c}_i \hat{v}_i \bigr\rangle \bigl\langle
X_n^c , \hat{c}_j \hat{v}_j
\bigr\rangle- \bigl\langle X_n^c , v_i \bigr
\rangle \bigl\langle X_n^c , v_j \bigr
\rangle \bigr)=\mathrm{O}_{P}\bigl(N^{-1/2}\bigr)
\end{equation}
and
%
\begin{equation}
\label{eqteaspoon} N^{-1/2} \sum_{n=1}^N
\varepsilon_n^{**} \bigl(\langle X_n - \bar
{X} , \hat {c}_i \hat{v}_i \rangle \langle
X_n - \bar{X} , \hat{c}_j \hat{v}_j
\rangle- \bigl\langle X_n^c , \hat{c}_i
\hat{v}_i \bigr\rangle \bigl\langle X_n^c ,
\hat{c}_j \hat {v}_j\bigr\rangle \bigr)=\mathrm{o}_{P}
(1 ).
\end{equation}
The left-hand side of \eqref{eqtablespoon} is
\[
N^{-1/2} \sum_{n=1}^N
\varepsilon_n^{**} \bigl\langle X_n^c
, \hat{c}_i \hat {v}_i \bigr\rangle \bigl(\bigl\langle
X_n^c , \hat{c}_j \hat{v}_j
\bigr\rangle- \bigl\langle X_n^c , v_j \bigr
\rangle \bigr) + N^{-1/2} \sum_{n=1}^N
\varepsilon_n^{**} \bigl\langle X_n^c
, v_j \bigr\rangle \bigl(\bigl\langle X_n^c
, \hat{c}_i \hat{v}_i \bigr\rangle- \bigl\langle
X_n^c , v_i \bigr\rangle \bigr).
\]
It follows from Assumptions \ref{agaussian}--\ref{aindep} that both
sets of random functions $\{ \varepsilon_nX_n^c(t)X_n^c(s), 1\leq n
\leq N \}$ and $\{ X_n^c(u)X_n^c(t)X_n^c(s), 1\leq n \leq N \}$ are
independent and identically distributed with zero mean so by the
central limit theorem in Hilbert spaces, we have
%
\begin{eqnarray}
\label{eqbosq1}
\Biggl\Vert N^{-1/2} \sum_{n=1}^N
\varepsilon_n X_n^c(t)X_n^c(s)
\Biggr\Vert_2&=&\mathrm{O}_{P}(1) \quad\mbox{and}\nonumber\\[-8pt]\\[-8pt]
\Biggl\Vert N^{-1/2} \sum
_{n=1}^N X_n^c(u)
X_n^c(t)X_n^c(s) \Biggr\Vert_2&=&\mathrm{O}_{P}(1).\nonumber
\end{eqnarray}
Next, we write that
\[
N^{-1/2} \sum_{n=1}^N
\varepsilon_n^{**} \bigl\langle X_n^c
, \hat{c}_i \hat {v}_i \bigr\rangle \bigl(\bigl\langle
X_n^c , \hat{c}_j \hat{v}_j
\bigr\rangle- \bigl\langle X_n^c , v_j \bigr
\rangle \bigr) = \delta_1 + \delta_2 +
\delta_3 + \delta_4,
\]
where, by \eqref{eqbosq1}, Lemma \ref{lgoodapproximation} and
repeated applications of the Cauchy--Schwarz inequality, we have
\begin{eqnarray*}
\llvert \delta_1\rrvert &=& \Biggl\llvert N^{-1/2} \sum
_{n=1}^N \varepsilon_n \bigl
\langle X_n^c , \hat{c}_i
\hat{v}_i \bigr\rangle \bigl(\bigl\langle X_n^c
, \hat{c}_j \hat{v}_j \bigr\rangle- \bigl\langle
X_n^c , v_j \bigr\rangle \bigr)\Biggr
\rrvert
\\
&\le&\Biggl\Vert  N^{-1/2} \sum
_{n=1}^N \varepsilon_n
X_n^c(t)X_n^c(s)
\hat{c}_i \hat{v}_i(t) \bigl( \hat{c}_j
\hat{v}_j(s) - v_j(s) \bigr) \Biggr
\Vert_1
\\
&\le&\Biggl\Vert  N^{-1/2} \sum
_{n=1}^N \varepsilon_n
X_n^c(t)X_n^c(s)
\Biggr\Vert_2 \bigl\Vert \hat{c}_j \hat
{v}_j(s) - v_j(s) \bigr\Vert_2
\\
&=&\mathrm{O}_{P}\bigl(N^{-1/2}\bigr),
\\
\llvert \delta_2\rrvert &=& \Biggl\llvert N^{-1/2} \sum
_{n=1}^N \sum_{k=p+1}^{\infty}
b_{k} \bigl\langle X_n^c, v_k
\bigr\rangle \bigl\langle X_n^c , \hat{c}_i
\hat {v}_i \bigr\rangle \bigl(\bigl\langle X_n^c
, \hat{c}_j \hat{v}_j \bigr\rangle- \bigl\langle
X_n^c , v_j \bigr\rangle \bigr)\Biggr
\rrvert
\\
&\le& \Biggl\Vert N^{-1/2} \sum
_{n=1}^N X_n^c(u)
X_n^c(t)X_n^c(s)
\Biggr\Vert_2 \Biggl\Vert  \sum
_{k=p+1}^{\infty} b_{k} v_k(u)\Biggr
\Vert _2 \bigl\Vert
\hat{c}_j \hat{v}_j(s) - v_j(s) \bigr\Vert_2
\\
&=&\mathrm{O}_{P}\bigl(N^{-1/2}\bigr),
\\
\llvert \delta_3\rrvert &=& \Biggl\llvert N^{-1/2} \sum
_{n=1}^N \sum_{k=1}^p
b_k \bigl\langle X_n^c, v_k -
\hat{c}_k \hat{v}_k \bigr\rangle \bigl\langle
X_n^c , \hat {c}_i \hat {v}_i
\bigr\rangle \bigl(\bigl\langle X_n^c ,
\hat{c}_j \hat{v}_j \bigr\rangle- \bigl\langle
X_n^c , v_j \bigr\rangle \bigr)\Biggr
\rrvert
\\
&\le&\sum_{k=1}^p \llvert
b_k\rrvert \Biggl\Vert  N^{-1/2} \sum
_{n=1}^N X_n^c(t)X_n^c(s)X_n^c(w)
\Biggr\Vert_2 \bigl\Vert
v_k(w) - \hat {c}_k \hat{v}_k(w) \bigr
\Vert _2 \bigl\Vert \hat{c}_j
\hat{v}_j(s) - v_j(s) \bigr\Vert_2
\\
&=&\mathrm{O}_{P}\bigl(N^{-1}\bigr),
\end{eqnarray*}
and
\begin{eqnarray*}
\llvert \delta_4\rrvert &=& \Biggl\llvert N^{-1/2} \sum
_{n=1}^N \sum_{k=1}^{p}
b_{k} \langle\bar{X} - \mu_X, \hat{c}_k
\hat{v}_k \rangle \bigl\langle X_n^c , \hat
{c}_i \hat{v}_i \bigr\rangle \bigl(\bigl\langle
X_n^c , \hat{c}_j \hat{v}_j
\bigr\rangle- \bigl\langle X_n^c , v_j \bigr
\rangle \bigr)\Biggr\rrvert
\\
&\le&\Biggl\llvert \sum_{k=1}^{p}
b_{k} \langle\bar{X} - \mu_X, \hat{c}_k \hat
{v}_k \rangle\Biggr\rrvert  \Biggl\Vert
N^{-1/2} \sum_{n=1}^N
X_n^c(t)X_n^c(s)\Biggr\Vert_2 \bigl\Vert \hat{c}_j \hat
{v}_j(s) - v_j(s) \bigr\Vert_2
\\
&=&\mathrm{O}_{P}\bigl(N^{-1/2}\bigr).
\end{eqnarray*}
Similarly,
\[
N^{-1/2} \sum_{n=1}^N
\varepsilon_n^{**} \bigl\langle X_n^c
, v_j \bigr\rangle \bigl(\bigl\langle X_n^c
, \hat{c}_i \hat{v}_i \bigr\rangle- \bigl\langle
X_n^c , v_i \bigr\rangle \bigr) =
\mathrm{o}_{P} (1 ),
\]
and therefore \eqref{eqtablespoon} is proven.

We now establish \eqref{eqteaspoon}.
The left-hand side of \eqref{eqteaspoon} is equal to
\[
N^{-1/2} \sum_{n=1}^N
\varepsilon_n^{**} \langle X_n - \bar{X} ,
\hat{c}_i \hat{v}_i \rangle \langle\mu_X -
\bar{X} , \hat{c}_j \hat{v}_j \rangle+ N^{-1/2}
\sum_{n=1}^N \varepsilon_n^{**}
\bigl\langle X_n^c , \hat{c}_j \hat
{v}_j \bigr\rangle \langle\mu_X - \bar{X} ,
\hat{c}_i \hat{v}_i \rangle.
\]
We write that
\[
N^{-1/2} \sum_{n=1}^N
\varepsilon_n^{**} \langle X_n - \bar{X} ,
\hat{c}_i \hat{v}_i \rangle \langle\mu_X -
\bar{X} , \hat{c}_j \hat{v}_j \rangle=
\delta_5 + \delta_6 + \delta_7 +
\delta_8,
\]
where, by the central limit theorem in Hilbert spaces, Lemma \ref
{lgoodapproximation}, and the Cauchy--Schwarz inequality, we have
\begin{eqnarray*}
\llvert \delta_5\rrvert &=&\Biggl\llvert N^{-1/2} \sum
_{n=1}^N \varepsilon_n \langle
X_n - \bar{X} , \hat{c}_i \hat{v}_i
\rangle \langle\mu_X - \bar {X} , \hat{c}_j
\hat{v}_j \rangle\Biggr\rrvert
\\
&\le&\bigl\llvert \langle\mu_X - \bar{X} ,
\hat{c}_j \hat{v}_j \rangle \bigr\rrvert
\Biggl\Vert N^{-1/2} \sum_{n=1}^N
\varepsilon_n \bigl(X_n(s) - \bar{X}(s) \bigr)
 \Biggr\Vert_2
\\
& =& \mathrm{O}_{P}\bigl(N^{-1/2}\bigr),
\\
\llvert \delta_6\rrvert &=&\Biggl\llvert N^{-1/2} \sum
_{n=1}^N \sum_{k=p+1}^{\infty}
b_{k} \bigl\langle X_n^c, v_k
\bigr\rangle \langle X_n - \bar {X} , \hat {c}_i
\hat{v}_i \rangle \langle\mu_X - \bar{X} ,
\hat{c}_j \hat {v}_j \rangle \Biggr\rrvert
\\
&\le&\bigl\llvert \langle\mu_X - \bar{X} , \hat{c}_j
\hat{v}_j \rangle \bigr\rrvert \Biggl\llvert N^{-1/2} \sum
_{n=1}^N \sum_{k=p+1}^{\infty}
b_{k} \bigl\langle X_n^c, v_k
\bigr\rangle \langle X_n - \bar{X} , \hat{c}_i
\hat{v}_i \rangle\Biggr\rrvert
\\
&=&\bigl\llvert \langle\mu_X - \bar{X} , \hat{c}_j
\hat{v}_j \rangle \bigr\rrvert \Biggl\llvert N^{-1/2} \sum
_{n=1}^N \Int\Int X_n^c(t)
\bigl( X_n(s) - \bar {X}(s) \bigr)\hat{v}_i(s) \sum
_{k=p+1}^{\infty} b_{k}v_k(t)
\,\mathrm{d}s\,\mathrm{d}t \Biggr\rrvert
\\
&=&\bigl\llvert \langle\mu_X - \bar{X} , \hat{c}_j
\hat{v}_j \rangle \bigr\rrvert \Biggl\llvert N^{-1/2} \sum
_{n=1}^N \Int\Int \bigl(X_n(t)-
\bar{X}(t) \bigr) \bigl( X_n(s) - \bar{X}(s) \bigr)\\
&&\hspace*{144pt}{}\times
\hat{v}_i(s) \sum_{k=p+1}^{\infty}
b_{k}v_k(t) \,\mathrm{d}s\,\mathrm{d}t\Biggr\rrvert
\\
&=& N^{1/2} \bigl\llvert \langle\mu_X - \bar{X} ,
\hat{c}_j \hat{v}_j \rangle \bigr\rrvert \Biggl\llvert
\Int\Int\hat{c}(t,s) \hat{v}_i(s) \sum_{k=p+1}^{\infty}
b_{k}v_k(t) \,\mathrm{d}s\,\mathrm{d}t\Biggr\rrvert
\\
&=& N^{1/2} \hat{\lambda}_i \bigl\llvert \langle
\mu_X - \bar{X} , \hat{c}_j \hat {v}_j
\rangle\bigr\rrvert \Biggl\llvert \Int\hat{v}_i(t) \sum
_{k=p+1}^{\infty} b_{k}v_k(t) \,\mathrm{d}t
\Biggr\rrvert
\\
&=& N^{1/2} \hat{\lambda}_i\bigl\llvert \langle
\mu_X - \bar{X} , \hat{c}_j \hat {v}_j
\rangle\bigr\rrvert \Biggl\llvert \Int\sum_{k=p+1}^{\infty}b_{k}
v_k(t) \bigl(\hat{v}_i(t) - \hat{c}_i
v_i(t) \bigr) \,\mathrm{d}t\Biggr\rrvert
\\
&\le& N^{1/2} \hat{\lambda}_i \bigl\llvert \langle
\mu_X - \bar{X} , \hat {c}_j \hat {v}_j
\rangle\bigr\rrvert  \Biggl\Vert \sum_{k=p+1}^{\infty}b_{k}
v_k(t) \Biggr\Vert_2 \bigl\Vert \hat{v}_i(t) - \hat{c}_i v_i(t)
 \bigr\Vert_2
\\
&=&\mathrm{O}_{P}\bigl(N^{-1/2}\bigr),
\\
\llvert \delta_7\rrvert &=&\Biggl\llvert N^{-1/2} \sum
_{n=1}^N \sum_{k=1}^p
b_k \bigl\langle X_n^c, v_k -
\hat{c}_k \hat{v}_k \bigr\rangle \langle
X_n - \bar {X} , \hat{c}_i \hat{v}_i
\rangle \langle\mu_X - \bar{X} , \hat{c}_j
\hat{v}_j \rangle\Biggr\rrvert
\\
&\le&\bigl\llvert \langle\mu_X - \bar{X} ,
\hat{c}_j \hat{v}_j \rangle \bigr\rrvert \Biggl\Vert
 N^{-1/2} \sum_{n=1}^N
\sum_{k=1}^p b_k
X_n^c(t) \bigl(X_n(s) - \bar{X}(s) \bigr)
\Biggr\Vert_2 \bigl\Vert
v_k(t) - \hat {c}_k \hat {v}_k(t)\bigr
\Vert_2
\\
& =& \mathrm{O}_{P}\bigl(N^{-1/2}\bigr)
\end{eqnarray*}
and
\begin{eqnarray*}
\llvert \delta_8\rrvert &=&\Biggl\llvert N^{-1/2} \sum
_{n=1}^N \sum_{k=1}^{p}
b_{k} \langle\bar{X} - \mu_X, \hat{c}_k
\hat{v}_k \rangle \langle X_n - \bar {X} ,
\hat{c}_i \hat{v}_i \rangle \langle\mu_X -
\bar{X} , \hat{c}_j \hat{v}_j \rangle \Biggr\rrvert
\\
&\le&\bigl\llvert \langle\mu_X - \bar{X} ,
\hat{c}_j \hat{v}_j \rangle \bigr\rrvert \Biggl\llvert
\sum_{k=1}^{p} b_{k} \langle
\bar{X} - \mu_X, \hat{c}_k \hat{v}_k
\rangle \Biggr\rrvert \Biggl\Vert  N^{-1/2} \sum
_{n=1}^N \bigl(X_n(s) - \bar{X}(s)
\bigr)\Biggr\Vert_2
\\
& =& \mathrm{O}_{P}\bigl(N^{-1/2}\bigr).
\end{eqnarray*}
This proves \eqref{eqteaspoon}, which also completes the proof of
\eqref{eqrunner} and hence \eqref{eqtofu}.

We proceed to the second step, which is the proof of \eqref
{eqicecream} and \eqref{eqorangejuice}. We express \eqref
{eqicecream} elementwise as
%
\begin{equation}
\label{eqcyclist} N^{-1/2} \sum_{n=1}^N
\bigl( \bigl\langle X_n^c, v_i \bigr\rangle
\bigl\langle X_n^c, v_j \bigr\rangle-
\lambda_i 1\{i=j\} \bigr) \Biggl( \sum_{k=1}^p
b_k \bigl\langle X_n^c, v_k -
\hat {c}_k \hat{v}_k \bigr\rangle \Biggr) =
\mathrm{o}_{P} (1 ).
\end{equation}
We observe that by the central limit theorem in Hilbert spaces and
Lemma \ref{lgoodapproximation} we have
\begin{eqnarray*}
\Biggl\llvert N^{-1/2} \sum_{n=1}^N
\Biggl( \sum_{k=1}^p b_k
\bigl\langle X_n^c, v_k -
\hat{c}_k \hat{v}_k \bigr\rangle \Biggr) \Biggr\rrvert
& \leq&\Biggl\Vert N^{-1/2} \sum
_{n=1}^N X_n^c(t)
\Biggr\Vert_2 \sum_{k=1}^p
|b_k| \bigl\Vert v_k(t) -
\hat{c}_k \hat{v}_k(t)  \bigr\Vert_2
\\
&=&\mathrm{O}_{P}\bigl(N^{-1/2}\bigr).
\end{eqnarray*}
Similarly,
\begin{eqnarray*}
&&\Biggl\llvert N^{-1/2} \sum_{n=1}^N
\bigl\langle X_n^c, v_i \bigr\rangle \bigl
\langle X_n^c, v_j \bigr\rangle \Biggl( \sum
_{k=1}^p b_k \bigl\langle
X_n^c, v_k - \hat{c}_k
\hat{v}_k \bigr\rangle \Biggr)\Biggr\rrvert
\\
&&\quad \le\sum_{k=1}^p \llvert
b_k\rrvert \Biggl\Vert n^{-1/2} \sum
_{n=1}^N X_n^c(t)X_n^c(s)X_n^c(w)
\Biggr\Vert_2 \bigl\Vert
v_k(w) - \hat{c}_k \hat{v}_k(w)\bigr
\Vert_2
\\
&&\quad=\mathrm{O}_{P}\bigl(N^{-1/2}\bigr).
\end{eqnarray*}
This proves \eqref{eqcyclist} and hence \eqref{eqicecream}. Next, we
establish \eqref{eqorangejuice}. We can express \eqref{eqorangejuice}
elementwise as
%
\begin{equation}
\label{eqscrewdriver} N^{-1/2} \sum_{n=1}^N
\bigl( \bigl\langle X_n^c, v_k \bigr\rangle
\bigl\langle X_n^c, v_{\ell} \bigr\rangle-
\lambda_k 1\{k={\ell}\} \bigr) \langle\bar{X} - \mu_X ,
\hat{c}_i \hat{v}_i \rangle=\mathrm{o}_{P} (1 ).
\end{equation}
Using the previous arguments, one can easily verify \eqref
{eqscrewdriver}, establishing \eqref{eqorangejuice}.

We will now finish the proof of the lemma by establishing \eqref
{eqredpen} as the third step. Using Assumptions \ref{agaussian}, \ref
{avarepsilon}, and \eqref{aindep}, we see that $\varepsilon_n^*$ has
mean zero and variance given by
\begin{eqnarray*}
E \bigl(\varepsilon_n^* \bigr)^2 &=& E \bigl(
\varepsilon_1^2 \bigr) + E \Biggl( \sum
_{i=p+1}^{\infty} \sum_{j=p+1}^{\infty}b_i
b_j \bigl\langle X^c_n, v_i
\bigr\rangle \bigl\langle X^c_n, v_j \bigr
\rangle \Biggr)
\\
&=& \sigma^2 + \sum_{i=p+1}^{\infty}b_i^2
E \bigl( \bigl\langle X^c_n, v_i \bigr
\rangle^2 \bigr)
\\
&=& \sigma^2 + \sum_{i=p+1}^{\infty}b_i^2
\lambda_i
\\
&=&\tau^2.
\end{eqnarray*}

Therefore, $\varepsilon_n^{*}  ({\mathbf D}_n - {\mathbf M} )$
is an iid sequence with mean zero and variance $\tau^2  ({\mathbf
G}-{\mathbf M}{\mathbf M^T} )$. The central limit theorem now
proves \eqref{eqredpen}, completing the proof of the lemma.
\end{pf}

\begin{lemma}\label{lestimatorsconverging}
If Assumptions \ref{aXmoments}--\ref{aunique} are satisfied, then
%
\begin{equation}
\label{eqsecret} %
\pmatrix{ \hat{\mathbf A}
\cr
\hat{\mathbf B}
\cr
\hat{\mu}}- %
\pmatrix{ {\mathbf0}
\cr
\tilde{\mathbf B}
\cr
\mu} =
\mathrm{O}_{P}\bigl(N^{-1/2}\bigr).
\end{equation}
In particular, we have
%
\begin{equation}
\label{eqbees} \bigl\| b_k v_k(t) - \hat{b}_k
\hat{v}_k(t) \bigr\|_2 = \mathrm{O}_{P}\bigl(N^{-1/2}
\bigr)
\end{equation}
and
%
\begin{equation}
\label{eqaees}\bigl \| \hat{a}_{i,j} \hat{v}_i(t)
\hat{v}_j(s) \bigr\|_2 = \mathrm{O}_{P}\bigl(N^{-1/2}
\bigr),
\end{equation}
where $\hat{a}_{i,j}$ and $\hat{b}_i$ are defined by
\[
\hat{\mathbf A}=\vech \bigl( \bigl\{ \hat{a}_{i,j}
\bigl(2-1\{i=j\} \bigr), 1\le i \le j\le p \bigr\}^T \bigr) \quad\mbox{and}\quad
\hat{\mathbf B}= ( \hat{b}_1, \hat{b}_2, \dots,
\hat{b}_p )^T.
\]
\end{lemma}

\begin{pf}
Lemmas \ref{lestimatingG} and \ref{lthelimitingdistribution} imply
that $\hat{\mathbf A} = \mathrm{O}_{P}(N^{-1/2})$. According to \eqref
{eqthehatsareback} and \eqref{eqrowsonly}, we can prove that
%
\begin{equation}
\label{eqthebeeestimate} \hat{\mathbf B}-\tilde{\mathbf B} = \mathrm{O}_{P}
\bigl(N^{-1/2}\bigr),
\end{equation}
by showing that
\[
\frac{1}{N} \sum_{n=1}^N
\varepsilon_n^{**} {\boldsymbol\Lambda}^{-1}
\hat{\mathbf F}_n = \mathrm{O}_{P}\bigl(N^{-1/2}\bigr)
\\
\]
or equivalently that
%
\begin{equation}
\label{eqbeeequive} \frac{1}{N} \sum_{n=1}^N
\varepsilon_n^{**} \langle X_n - \bar{X},
\hat {v}_i \rangle
\\
= \mathrm{O}_{P}\bigl(N^{-1/2}\bigr).
\\
\end{equation}
We note that
\[
\frac{1}{N} \sum_{n=1}^N
\varepsilon_n^{**} \langle X_n - \bar{X}, \hat
{v}_i \rangle
\\
= \delta_9 + \delta_{10} + \delta_{11} +
\delta_{12},
\\
\]
where, following the arguments in the proof of Lemma \ref
{lthelimitingdistribution}, one can verify that
\begin{eqnarray*}
\llvert \delta_{9}\rrvert &=&\Biggl\llvert \frac{1}{N} \sum
_{n=1}^N \varepsilon_n \langle
X_n - \bar{X}, \hat{v}_i \rangle\Biggr\rrvert
\mathrm{O}_{P}\bigl(N^{-1/2}\bigr),
\\
\llvert \delta_{10}\rrvert &=&\Biggl\llvert \frac{1}{N} \sum
_{n=1}^N \sum_{k=p+1}^{\infty}
b_{k} \bigl\langle X_n^c, v_k
\bigr\rangle \langle X_n - \bar {X}, \hat{v}_i \rangle
\Biggr\rrvert =\mathrm{O}_{P}\bigl(N^{-1/2}\bigr),
\\
\llvert \delta_{11}\rrvert &=&\Biggl\llvert \frac{1}{N} \sum
_{n=1}^N \sum_{k=1}^p
b_k \bigl\langle X_n^c, v_k -
\hat{c}_k \hat{v}_k \bigr\rangle \langle
X_n - \bar{X}, \hat {v}_i \rangle\Biggr\rrvert
=\mathrm{O}_{P}\bigl(N^{-1/2}\bigr)
\end{eqnarray*}
and
\[
\llvert \delta_{12}\rrvert =\Biggl\llvert \frac{1}{N} \sum
_{n=1}^N \sum_{k=1}^{p}
b_{k} \langle\bar{X} - \mu_X, \hat{c}_k
\hat{v}_k \rangle \langle X_n - \bar{X},
\hat{v}_i \rangle\Biggr\rrvert %
=\mathrm{O}_{P}\bigl(N^{-1/2}\bigr). %
\]

This proves \eqref{eqbeeequive} and hence \eqref{eqthebeeestimate}.

To complete the justification of \eqref{eqsecret}, we need to show that
%
\begin{equation}
\label{eqmuestimate} \hat{\mu}- \mu= \mathrm{O}_{P}\bigl(N^{-1/2}
\bigr).
\end{equation}
Due to \eqref{eqthehatsareback} and \eqref{eqrowsonly}, \eqref
{eqmuestimate} will be established by proving that
%
\begin{equation}
\label{eqthefirstnoel} \frac{1}{N}\sum_{n=1}^{N}
\varepsilon_n^{**} \bigl(-{\mathbf M}^T \bigl({
\mathbf G}-{\mathbf M} {\mathbf M}^T \bigr)^{-1}\hat{\mathbf
D}_n + 1+{\mathbf M}^T \bigl({\mathbf G}-{\mathbf M} {
\mathbf M}^T \bigr)^{-1} {\mathbf M} \bigr) =
\mathrm{O}_{P}\bigl(N^{-1/2}\bigr).
\end{equation}

To prove \eqref{eqthefirstnoel}, it is sufficient to show
%
\begin{equation}
\label{eqthesecondnoel} \frac{1}{N}\sum_{n=1}^{N}
\varepsilon_n^{**} \hat{\mathbf D}_n =
\mathrm{O}_{P}\bigl(N^{-1/2}\bigr)
\end{equation}
and
%
\begin{equation}
\label{eqthethirdnoel} \frac{1}{N}\sum_{n=1}^{N}
\varepsilon_n^{**}= \mathrm{O}_{P}\bigl(N^{-1/2}
\bigr).
\end{equation}
Due to Lemma \ref{lthelimitingdistribution}, \eqref{eqthethirdnoel}
implies \eqref{eqthesecondnoel}, so we prove only \eqref
{eqthethirdnoel}. We write that
\[
\frac{1}{N}\sum_{n=1}^{N}
\varepsilon_n^{**} = \delta_{13}+
\delta_{14}+\delta_{15}+\delta_{16},
\]
where, by the central limit theorem in Hilbert spaces and Lemma \ref
{lgoodapproximation}, we have
\begin{eqnarray*}
\llvert \delta_{13}\rrvert &=&\Biggl\llvert \frac{1}{N}\sum
_{n=1}^{N} \varepsilon_n\Biggr
\rrvert %
=\mathrm{O}_{P}\bigl(N^{-1/2}\bigr),
\\
\llvert \delta_{14}\rrvert &=&\Biggl\llvert
\frac{1}{N}\sum_{n=1}^{N} \sum
_{k=p+1}^{\infty} b_{k} \bigl\langle
X_n^c, v_k \bigr\rangle\Biggr\rrvert
\le\Biggl\Vert  \frac{1}{N}\sum
_{n=1}^{N} X_n^c(t)\Biggr
\Vert_2\Biggl\Vert  \sum
_{k=p+1}^{\infty} b_{k} v_k(t)
\Biggr\Vert_2 %
=\mathrm{O}_{P}
\bigl(N^{-1/2}\bigr),
\\
\llvert \delta_{15}\rrvert &=&\Biggl\llvert
\frac{1}{N}\sum_{n=1}^{N} \sum
_{k=1}^p b_k \bigl\langle
X_n^c, v_k - \hat{c}_k
\hat{v}_k(t) \bigr\rangle\Biggr\rrvert %
=\mathrm{O}_{P}\bigl(N^{-1}\bigr)
\end{eqnarray*}
and
\[
\llvert \delta_{16}\rrvert =\Biggl\llvert \frac{1}{N}\sum
_{n=1}^{N} \sum
_{k=1}^{p} b_{k} \langle\bar{X} -
\mu_X, \hat{c}_k \hat{v}_k \rangle\Biggr
\rrvert
=\mathrm{O}_{P}\bigl(N^{-1/2}\bigr).
\]
This proves \eqref{eqthethirdnoel}, which establishes \eqref
{eqmuestimate} and completes the proof of \eqref{eqsecret}.

Using \eqref{eqsecret} and Lemma \ref{lgoodapproximation}, we will
now show \eqref{eqbees} and \eqref{eqaees}. We conclude from \eqref
{eqsecret} that
\[
\hat{b}_i - \hat{c}_i b_i =
\mathrm{O}_{P}\bigl(N^{-1/2}\bigr)\quad \mbox{and} \quad\hat{a}_{i,j} =
\mathrm{O}_{P}\bigl(N^{-1/2}\bigr).
\]
Now, Lemma \ref{lgoodapproximation} yields that
\begin{eqnarray*}
\bigl\| b_k v_k(t) - \hat{b}_k
\hat{v}_k(t) \bigr\|_2 &\le&\bigl\| b_k
\bigl(v_k(t) -\hat {c}_k \hat{v}_k(t)\bigr)
\bigr\|_2 + \bigl\|(b_k\hat{c}_k -
\hat{b}_k )\hat {v}_k(t) \bigr\|_2
\\
&\le&|b_k| \bigl\| v_k(t) -\hat{c}_k
\hat{v}_k(t) \bigr\|_2 + |b_k
\hat{c}_k - \hat {b}_k|
\\
&=& \mathrm{O}_{P}\bigl(N^{-1/2}\bigr).
\end{eqnarray*}
Similarly,
\[
\bigl\| \hat{a}_{i,j} \hat{v}_i(t)
\hat{v}_j(s)\bigr \|_2 %
= \mathrm{O}_{P}
\bigl(N^{-1/2}\bigr). %
\]
This proves \eqref{eqbees} and \eqref{eqaees} and completes the proof
of the lemma.
\end{pf}

\begin{lemma}\label{ltau}
If Assumptions \ref{agaussian}--\ref{aunique} are satisfied, then
\[
\hat{\tau}^2 - \tau^2 = \mathrm{O}_{P}
\bigl(N^{-1/2}\bigr).
\]
\end{lemma}
\begin{pf}
Since
\[
\frac{1}{N} \sum_{n=1}^{N}
\varepsilon_n^{*2} - \tau^2 \as0,
\]
it is enough to show that
%
\begin{equation}
\label{eqsquaresareclose} \frac{1}{N} \sum_{n=1}^{N}
\bigl(\hat{\varepsilon}_n^2 - \varepsilon_n^{*2}
\bigr) = \mathrm{O}_{P}\bigl(N^{-1/2}\bigr).
\end{equation}
Since
\[
\frac{1}{N} \sum_{n=1}^{N} \bigl(
\hat{\varepsilon}_n^2 - \varepsilon_n^{*2}
\bigr) = \frac{1}{N} \sum_{n=1}^{N}
\bigl(\hat{\varepsilon }_n - \varepsilon_n^{*}
\bigr) \bigl(\hat{\varepsilon}_n + \varepsilon_n^{*}
\bigr)%
=\frac{1}{N} \sum_{n=1}^{N}
\bigl(\hat{\varepsilon}_n - \varepsilon_n^{*}
\bigr) \hat{\varepsilon}_n + \frac{1}{N} \sum
_{n=1}^{N} \bigl(\hat{\varepsilon}_n -
\varepsilon_n^{*} \bigr) \varepsilon_n^{*},
\]
\eqref{eqsquaresareclose} follows from
%
\begin{equation}
\label{eqleftforthereader} \Biggl\llvert \frac{1}{N} \sum
_{n=1}^{N} \bigl(\hat{\varepsilon}_n -
\varepsilon_n^{*} \bigr) \varepsilon_n^{*}
\Biggr\rrvert = \mathrm{O}_{P}\bigl(N^{-1/2}\bigr)
\end{equation}
and
%
\begin{equation}
\label{eqIdoit} \Biggl\llvert \frac{1}{N} \sum
_{n=1}^{N} \bigl(\hat{\varepsilon}_n -
\varepsilon_n^{*} \bigr) \hat{\varepsilon}_n
\Biggr\rrvert = \mathrm{O}_{P}\bigl(N^{-1/2}\bigr).
\end{equation}

We decompose \eqref{eqleftforthereader} as
\[
\frac{1}{N} \sum_{n=1}^{N} \bigl(
\hat{\varepsilon}_n - \varepsilon_n^{*}
\bigr) {\varepsilon}_n^* = \eta_1 + \eta_2 +
\eta_3,
\]
where
\begin{eqnarray*}
\eta_1 &=&\frac{1}{N} \sum_{n=1}^N
\varepsilon_n^* (\mu- \hat {\mu } ),
\\
\eta_2 &=&\frac{1}{N} \sum_{n=1}^N
\varepsilon_n^*\sum_{i=1}^p
\bigl( b_i \bigl\langle X_n^c,
v_i \bigr\rangle- \hat{b}_i \langle X_n -
\bar{X}, \hat{v}_i\rangle \bigr),
\\
\eta_3 &=&\frac{1}{N} \sum_{n=1}^N
\varepsilon_n^*\sum_{i=1}^p
\sum_{j=i}^p \bigl(2-1\{i=j\}\bigr) \bigl(
a_{i,j}\bigl\langle X_n^c, v_i
\bigr\rangle \bigl\langle X_n^c, v_j \bigr
\rangle - \hat{a}_{i,j} \langle X_n -\bar{X} ,
\hat{v}_i\rangle \langle X_n -\bar{X} , \hat
{v}_j\rangle \bigr).
\end{eqnarray*}
It is clear that $\eta_1=\mathrm{O}_{P}(N^{-1})$. We also see that $\eta_2 =
\eta_{2,1} + \eta_{2,2} + \eta_{2,3} + \eta_{2,4}$,
where
\begin{eqnarray*}
\eta_{2,1}&=&\frac{1}{N} \sum_{n=1}^N
Y_n \sum_{i=1}^p
\bigl(b_i\bigl\langle X_n^c,
v_i\bigr\rangle- \hat{b}_i \langle X_n -
\bar{X}, \hat{v}_i \rangle \bigr),
\\
\eta_{2,2}&=&-\frac{1}{N} \sum_{n=1}^N
\mu\sum_{i=1}^p \bigl(b_i
\bigl\langle X_n^c, v_i\bigr\rangle-
\hat{b}_i \langle X_n - \bar{X}, \hat{v}_i
\rangle \bigr),
\\
\eta_{2,3}&=&-\frac{1}{N} \sum_{n=1}^N
\sum_{\ell=1}^p b_{\ell} \bigl\langle
X_n^c, v_{\ell} \bigr\rangle\sum
_{i=1}^p \bigl(b_i\bigl\langle
X_n^c, v_i\bigr\rangle-
\hat{b}_i \langle X_n - \bar{X}, \hat{v}_i
\rangle \bigr),
\\
\eta_{2,4}&=&-\frac{1}{N} \sum_{n=1}^N
\sum_{\ell=1}^p \sum
_{k=\ell}^p \bigl(2-1\{k=\ell\}\bigr) a_{\ell, k}
\bigl\langle X_n^c, v_{\ell}\bigr\rangle \bigl
\langle X_n^c, v_{k}\bigr\rangle \sum
_{i=1}^p \bigl(b_i\bigl\langle
X_n^c, v_i\bigr\rangle-
\hat{b}_i \langle X_n - \bar{X}, \hat{v}_i
\rangle \bigr).
\end{eqnarray*}
Applying \eqref{eqbees} and the central limit theorem in Hilbert
spaces, we obtain that
\begin{eqnarray*}
\llvert \eta_{2,1}\rrvert &=&\Biggl\llvert \frac{1}{N} \sum
_{n=1}^N Y_n \sum
_{i=1}^p \bigl(b_i\bigl\langle
X_n^c, v_i\bigr\rangle-
\hat{b}_i \langle X_n - \bar{X}, \hat {v}_i
\rangle \bigr)\Biggr\rrvert
\\
&\le&\sum_{i=1}^p\Biggl\Vert  \frac{1}{N} \sum_{n=1}^N
Y_n \bigl(b_i X_n^c(t)
v_i(t) - \hat{b}_i \bigl(X_n(t) -
\bar{X}(t) \bigr) \hat {v}_i(t) \bigr)\Biggr\Vert_1
\\
&\le&\sum_{i=1}^p\Biggl\Vert  \frac{1}{N} \sum_{n=1}^N
Y_n X_n(t) \bigl(b_i v_i(t) -
\hat{b}_i \hat{v}_i(t) \bigr)\Biggr\Vert_1
\\
&&{}+ \sum_{i=1}^p\Biggl\Vert \frac{1}{N} \sum_{n=1}^N
Y_n \bigl(b_i \mu_X(t) v_i(t)
- \hat{b}_i \bar{X}(t) \hat{v}_i(t) \bigr)
\Biggr\Vert_1
\\
&\le&\sum_{i=1}^p \Biggl
\Vert \frac{1}{N} \sum_{n=1}^N
Y_n X_n(t) \Biggr\Vert_2
\bigl\Vert b_i v_i(t) -
\hat{b}_i \hat{v}_i(t)  \bigr\Vert_2
\\
&&{}+ \sum_{i=1}^p\Biggl
\Vert \frac{1}{N} \sum_{n=1}^N
Y_n \bar {X}(t) \bigl(\hat{b}_i \hat{v}_i(t)
- b_i v_i(t) \bigr)\Biggr
\Vert_1
\\
&&{}+ \sum_{i=1}^p \Biggl
\Vert \frac{1}{N} \sum_{n=1}^N
Y_n b_i v_i(t) \bigl(\bar{X}(t) -
\mu_X(t) \bigr)\Biggr\Vert_1
\\
&\le&\sum_{i=1}^p\Biggl
\Vert \frac{1}{N} \sum_{n=1}^N
Y_n X_n(t)\Biggr\Vert_2
\bigl\Vert b_i v_i(t) -
\hat{b}_i \hat{v}_i(t)  \bigr\Vert_2
\\
&&{}+ \sum_{i=1}^p \Biggl
\Vert \frac{1}{N} \sum_{n=1}^N
Y_n \bar{X}(t)\Biggr\Vert_2\bigl\Vert \hat{b}_i \hat{v}_i(t) - b_i
v_i(t)  \bigr\Vert_2
\\
&&{}+ \sum_{i=1}^p \Biggl
\Vert \frac{1}{N} \sum_{n=1}^N
Y_n b_i v_i(t) \Biggr
\Vert_2 \bigl\Vert \bar{X}(t) - \mu(t) \bigr\Vert_2
\\
&=&\mathrm{O}_{P}\bigl(N^{-1/2}\bigr).
\end{eqnarray*}
In a like manner, one can verify that $\eta_{2,i} = \mathrm{O}_{P}(N^{-1/2}), i=2,3,4.$

This proves that $\eta_2=\mathrm{O}_{P}(N^{-1/2})$. In a similar fashion, one can
show that $\eta_3 =\mathrm{O}_{P}(N^{-1/2})$. This proves \eqref
{eqleftforthereader}. Following the previous arguments, one can
establish \eqref{eqIdoit}, completing the proof of the lemma.
\end{pf}

\section*{Acknowledgement}
Supported in part by NSF Grant DMS-09-05400.



\printhistory

\end{document}